\documentstyle[12pt,amssymb,amsbsy,righttag]{amsart}
\renewcommand{\a}{\alpha}

\newcommand{\g}{\gamma}
\renewcommand{\d}{\delta}

\newcommand{\z}{\zeta}

\newcommand{\r}{\rho}
\newcommand{\s}{\sigma}

\newcommand{\f}{\varphi}

\newcommand{\D}{\Delta}
\renewcommand{\L}{\Lambda}

\renewcommand{\O}{\Omega}
\newcommand{\U}{\Upsilon}

\newcommand{\A}{{\cal A}}

\newcommand{\cL}{{\cal L}}
\newcommand{\M}{{\cal M}}

\newcommand{\cR}{{\cal R}}

\newcommand{\V}{{\cal V}}
\newcommand{\W}{{\cal W}}

\newcommand{\C}{{\Bbb C}}
\newcommand{\T}{{\Bbb T}}
\newcommand{\pp}{{\Bbb P}}
\newcommand{\dd}{{\Bbb D}}
\newcommand{\R}{{\Bbb R}}
\newcommand{\Z}{{\Bbb Z}}
\newcommand{\mm}{{\Bbb M}}
\newcommand{\0}{{\Bbb O}}

\newcommand{\PP}{{\boldsymbol{\cal P}}}

\newcommand{\m}{{\boldsymbol m}}

\newcommand{\rf}[1]{(\ref{#1})}

\newcommand{\df}{\stackrel{\mathrm{def}}{=}}
\newcommand{\dist}{\operatorname{dist}}
\newcommand{\Ker}{\operatorname{Ker}}

\newcommand{\clos}{\operatorname{clos}}

\newcommand{\rank}{\operatorname{rank}}
\newcommand{\const}{\operatorname{const}}

\newcommand{\eeq}{\end{equation}}
\newcommand{\beq}{\begin{equation}}
\newcommand{\bay}{\begin{eqnarray}}
\newcommand{\ey}{\end{eqnarray}}

\newcommand{\be}{\infty}

\newcommand{\bl}{\blacksquare}
\newcommand{\ess}{\operatorname{ess}}
\newcommand{\ind}{\operatorname{ind}}
\newcommand{\wind}{\operatorname{wind}}
\newcommand{\Range}{\operatorname{Range}}

\newcommand{\Pf}{{\bf Proof. }}

\newcommand{\ov}{\overline}

\newtheorem{thm}{\hspace{\parindent}Theorem}[section]

\newtheorem{cor}[thm]{\hspace{\parindent}Corollary}
\newtheorem{lem}[thm]{\hspace{\parindent}Lemma}

\begin{document}
\newcommand{\vse}{\vspace{.2in}}
\newcommand{\bs}{\boldsymbol}
\renewcommand{\theequation}{\thesection.\arabic{equation}}

\author{R.B. Alexeev and V.V. Peller}
\thanks{The second author is partially supported by NSF grant DMS 9970561}

\title{Badly approximable matrix functions \\and canonical factorizations}

\maketitle

\begin{abstract}
We continue studying the problem of analytic approximation of matrix 
functions. We introduce the notion of a partial canonical factorization
of a badly approximable matrix function $\Phi$ and the notion of a canonical 
factorization of a very badly approximable matrix function $\Phi$. 
Such factorizations  are defined in terms of so-called balanced unitary-valued 
functions which have many remarkable properties.
Unlike the case of thematic factorizations studied earlier in [PY1], [PY2], [PT], 
[AP1], the factors in
canonical factorizations (as well as partial canonical factorizations) are uniquely
determined by the matrix function $\Phi$ up to constant unitary factors.
We study many properties of canonical factorizations. In particular we show that
under certain natural assumptions on a function space $X$ the condition
$\Phi\in X$ implies that all factors in a canonical factorization of $\Phi$ 
belong to the same space $X$. In the last section we characterize the very badly
approximable unitary-valued functions $U$ that satisfy the condition 
$\|H_U\|_{\text e}<1$.
\end{abstract}

\setcounter{equation}{0}
\setcounter{section}{0}
\section{\bf Introduction}

\

The problem of uniform approximation by bounded analytic functions has been
studied for a long time. It was proved in [Kh] that for a continuous function $\f$
on the unit circle $\T$ there exists a unique best approximation $f$ by bounded
analytic functions and the error function $\f-f$ has constant modulus almost 
everywhere on $\T$. Since that time this approximation problem has been studied
by many authors. An important step in developing the theory of approximation by 
analytic functions was the Nehari theorem [Ne] according to which the distance from
an $L^\be$ function $\f$ to the space $H^\be$ of bounded analytic functions
is equal to the norm of the {\it Hankel operator} 
$H_\f:H^2\to H^2_-\df L^2\ominus H^2$ defined by
$$
H_\f f=\pp_-\f f,
$$
where $\pp_-$ is the orthogonal projection onto $H^2_-$. Hankel operators were 
used essentially for further development of this theory in [AAK1-3] and [PK].

Later it turned out that the approximation problem in question
plays a crucial role in so-called $H^\be$ control theory. Moreover, for the needs
of control theory engineers have to consider matrix-valued functions.
We refer the reader to [F] for an introduction in $H^\be$ control.

In this paper we continue the study of best analytic approximation of matrix-valued
functions. We consider the space $L^\be(\mm_{m,n})$ of essentially bounded functions
which take values in the space $\mm_{m,n}$ of $m\times n$ matrices and endow it with 
the norm
$$
\|\Phi\|_{L^\be}=\ess\sup_{\z\in\T}\|\Phi(\z)\|_{\mm_{m,n}},
$$ 
(the space $\mm_{m,n}$ is endowed with the operator norm in the space of operators 
from $\C^n$ to $\C^m$) and we study approximations of functions 
$\Phi\in L^\be(\mm_{m,n})$ by functions in the subspace $H^\be(\mm_{m,n})$ of 
$L^\be(\mm_{m,n})$ that consists of bounded analytic matrix functions in the unit 
disk $\dd$.

However, it is well known and it is easy to see that in the matrix case a continuous
(and even infinitely smooth) function $\Phi$ generically has infinitely many
best approximations by bounded analytic functions. It seems natural to impose 
additional constraints on a best approximation and choose among best approximations 
the ``very best''.

To introduce the notion of very best approximation, we recall that for a matrix $A$
(or a Hilbert space operator $A$) the singular value $s_j(A)$, $j\ge0$, is by 
definition the distance from $A$  to the set of matrices (operators) of rank at
most $j$. Clearly, $s_0(A)=\|A\|$.

Given a matrix function $\Phi\in L^\be(\mm_{m,n})$ we define inductively
the sets $\bs{\O}_j$, \linebreak$0\le j\le\min\{m,n\}-1$, by
$$
\bs{\O}_0=\{F\in H^\be(\mm_{m,n})
:~F~\mbox{minimizes}~\ t_0\df\ess\sup_{\z\in\T}\|\Phi(\z)-F(\z)\|\};
$$
$$
\bs{\O}_j=\{F\in \O_{j-1}:~F~\mbox{minimizes}~\ 
t_j\df\ess\sup_{\z\in\T}s_j(\Phi(\z)-F(\z))\}.
$$
Functions in $\bs{\O}_{\min\{m,n\}-1}$ are called {\it superoptimal approximations}
of $\Phi$ by bounded analytic matrix functions. The numbers $t_j=t_j(\Phi)$
are called the {\it superoptimal singular values} of $\Phi$. Note that the functions
in $\bs{\O}_0$ are just the best approximations by analytic matrix functions.
The notion of superoptimal approximation was introduced in [Y].

It was proved in [PY1] that for a continuous $m\times n$ matrix function $\Phi$
there exists a unique superoptimal approximation $F$ by bounded analytic matrix
functions and the error function $\Phi-F$ satisfies
\beq
\label{1.1}
s_j(\Phi(\z)-F(\z))=t_j(\Phi)\quad\mbox{almost everywhere on}\quad\T.
\end{equation}
Later this uniqueness result was obtained in a different way by Treil [T].

In [PT] the uniqueness result was improved. It was shown in [PT] that if
\linebreak$\Phi\in L^\be(\mm_{m,n})$ and the essential norm $\|H_\Phi\|_{\text e}$
of the Hankel operator $H_\Phi$ is less than the smallest nonzero superoptimal 
singular value of $\Phi$, then there exists a unique superoptimal approximation $F$ 
by analytic matrix functions and \rf{1.1} holds. 

For $\Phi\in L^\be(\mm_{m,n})$ the Hankel operator
$H_\Phi:H^2(\C^n)\to H^2_-(\C^n)$ is defined in the same way as in the scalar
case: $H_\Phi f=\pp_-\Phi f$, its norm is given by
$$
\|H_\Phi\|=\dist_{L^\be}(\Phi,H^\be(\mm_{m,n})),
$$
([Pa]) and its essential norm is equal to
$$
\|H_\Phi\|_{\text e}=\dist_{L^\be}(\Phi,(H^\be+C)(\mm_{m,n}))
$$
(see e.g., [Sa] where the proof of this formula in the scalar case is given, 
in the matrix case the proof is the same). Clearly, $\|H_\Phi\|$ is equal to 
the largest superoptimal singular value $t_0(\Phi)$ of $\Phi$.

In [PY1] and [PT] it was shown that if $\Phi$ satisfies the above conditions and
$F$ is the unique superoptimal approximation of $\Phi$,
then the error function $\Phi-F$ admits a so-called thematic factorization
(see \S 2 for precise definitions). The technique of thematic factorizations
turned out to be very fruitful (see [PY1], [PY2], [PT]). However, in the case
of multiple superoptimal singular values the factors in a thematic factorization
essentially depend on the choice of a factorization and by no means they can
be determined by the matrix function $\Phi$ itself.

In this paper we consider a modification of the notion of a thematic factorization.
Under the same assumptions on $\Phi$ we show that for the superoptimal approximation
$F$ the error function $\Phi-F$ admits a so-called canonical factorization. Unlike
the case of thematic factorizations the factors in canonical factorizations
are determined by the function $\Phi$ modulo constant unitary factors.
Similarly, we consider so-called partial canonical factorizations for badly
approximable matrix functions and prove the same invariance properties.
This is done in \S 8.

Canonical factorizations are defined in terms of so called balanced unitary-valued
matrix functions which are defined in \S 3. The notion of a balanced matrix function  
generalizes the notion of a thematic matrix function that was used
to define a thematic factorization. We discuss some remarkable properties of balanced 
matrix functions in \S 3.

Recall that a matrix function $\Phi$ is called {\it badly approximable} if the zero 
function is a best approximation of $\Phi$ by bounded analytic matrix functions.
A matrix function $\Phi$ is called {\it very badly approximable} if the zero 
function is a superoptimal approximation of $\Phi$ by bounded analytic matrix 
functions.

In \S 4 we prove that badly approximable matrix functions admit so-called
partial canonical factorizations. In \S 5 and \S 6 we compare (partial) thematic
factorizations with (partial) canonical factorizations and deduce a number of results
on partial canonical factorizations from the corresponding results on
partial thematic factorizations proved earlier in [PY1] and [PT].
Canonical factorizations of very badly approximable functions are discussed in
\S 7. 

\S 9 is devoted to hereditary properties of (partial) canonical
factorizations. In other words for many important function spaces $X$ we prove 
in \S 9 that if  $\Phi$ is a badly approximable matrix function whose
entries belong to $X$, then the entries of all factors in a partial canonical
factorization of $\Phi$ also belong to $X$. In particular if 
$\Phi\in(H^\be+C)(\mm_{m,n})$, then the entries of all factors in a canonical
factorization of $\Phi$ belong to the space $QC$ of quasi-continuous functions.

In \S 10 we characterize the very badly approximable unitary-valued
functions $U$ under the assumption $\|H_U\|_{\text e}<1$. Such unitary-valued
functions are involved in canonical factorizations.

Finally, in \S 2 we give definitions and state results that will be used in 
this paper.

\

\setcounter{equation}{0}
\setcounter{section}{1}
\section{\bf Preliminaries}

\

{\bf Toeplitz operators and Wiener--Hopf factorizations.} For a matrix function 
$\Phi\in L^\be(\mm_{m,n})$ the {\it Toeplitz operator} 
$T_\Phi:H^2(\C^n)\to H^2(\C^m)$ is defined by
$$
T_\Phi f=\pp_+\Phi f,\quad f\in H^2(\C^n).
$$
Suppose now that $m=n$. By Simonenko's theorem [Si] (see also [LS]), if $T_\Phi$ is 
Fredholm, then $\Phi$ admits a Wiener--Hopf factorization
$$
\Phi=Q^*_2\left(\begin{array}{cccc}z^{d_1}&0&\cdots&0\\
0&z^{d_2}&\cdots&0\\\vdots&\vdots&\ddots&\vdots\\0&0&\cdots&z^{d_n}
\end{array}\right)Q_1^{-1},
$$
where $d_1,\cdots,d_n\in\Z$, and $Q_1$ and $Q_2$ are matrix functions invertible
in $H^2(\mm_{n,n})$. It is always possible to arrange the {\it Wiener--Hopf indices}
$d_j$ is the nondecreasing order: $d_1\le\cdots\le d_n$ in which case they are 
uniquely determined by $\Phi$. 

{\bf Maximizing vectors of vectorial Hankel operators.} Let $\Phi$ be a matrix 
function in $L^\be(\mm_{m,n})$ such that the Hankel operator 
$H_\Phi:H^2(\C^n)\to H^2_-(\C^m)$ has a maximizing vector $f\in H^2(\C^n)$.
Let $F$ be a best approximation of $\Phi$ by bounded analytic matrix functions.
Put $g=H_\Phi f$. Then
$$
H_\Phi f=(\Phi-F)f,
$$
$$
\|\Phi(\z)-F(\z)\|_{\mm_{m,n}}=\|H_\Phi\|\quad\mbox{for almost all}\quad\z\in\T,
$$
$$
\|g(\z)\|_{\C^m}=\|H_\Phi\|\cdot\|f(\z)\|_{\C^n}\quad\mbox{for almost all}
\quad\z\in\T,
$$
and $f(\z)$ is a maximizing vector of $\Phi(\z)-F(\z)$ for almost all $\z\in\T$
([AAK3], see also [PY1]).

{\bf Badly approximable scalar functions.} If $\f$ is a 
nonzero continuous scalar function on the unit circle, then $\f$ is badly 
approximable if and only if $\f$ has constant modulus and its winding number 
$\wind\f$ is negative (see [AAK1] and [Po]). If $\f$ is a nonzero scalar function in 
$H^\be+C$, then $\f$ is badly approximable if and only if $|\f|$ is constant almost
everywhere on $\T$, $\f\in QC$ and $\ind T_\f>0$ (see [PK]). Here 
$$
H^\be+C\df\{f+g:~f\in C(\T),~g\in H^\be\}
$$
is a closed subalgebra of $L^\be$ and 
$$
QC\df\{f\in H^\be+C:~\bar f\in H^\be+C\}.
$$
Note that if $\f\in QC$ and $|\f|=\const>0$, then the Toeplitz operator
$T_\f$ is Fredholm.

These results can easily be generalized to the set of scalar functions $\f\in L^\be$ 
satisfying $\|H_\f\|_{\text e}<\|H_\f\|$. Under this condition $\f$ is
badly approximable if and only if $|\f|=\const$ and $\ind T_\f>0$. Again
it is easy to see that if $|\f|=\const$ and $\|H_\f\|_{\text e}<\|H_\f\|$, then
$T_\f$ is Fredholm.

{\bf Inner and outer matrix functions.} A matrix function $\Phi\in H^\be(\mm_{m,n})$
is called {\it inner} if $\Phi^*(\z)\Phi(\z)=I_n$ for almost all $\z\in\T$, where 
$I_n$ stands for the identity $n\times n$ matrix (or the matrix function
identically equal to $I_n$). 

Consider the operator of multiplication by $z$ on $H^2(\C^n)$. If $\cL$ is a nonzero
invariant subspace of this operator, then by the Beurling--Lax theorem, there exists 
an inner matrix function $\Theta\in H^\be(\mm_{n,r})$ such that
$$
\cL=\Theta H^2(\C^r).
$$
In this case
$$
\dim\{f(\z):~f\in\cL\}=r\quad\mbox{for almost all}\quad\z\in\dd.
$$
If $\Theta^\circ$ is an inner matrix function in $H^\be(\mm_{n,r^\circ})$ and
$\Theta^\circ H^2(\C^{r^\circ})=\Theta H^2(r)$, then $r=r^\circ$ and
there exists a unitary matrix $\frak U\in\mm_{r,r}$ such that 
$\Theta^\circ=\Theta\frak U$.

A matrix function $F\in H^2(\mm_{m,n})$ is called {\it outer} if the linear span
of the set $\{Fz^jx:~j\ge0,~x\in\C^n\}$ is dense in $H^2(\C^m)$. $F$ is called 
{\it co-outer} if the transposed function $F^{\text t}$ is outer.

If $\Psi$ is a matrix function in $H^2(\mm_{m,n})$, then there exist an inner matrix 
function $\Theta$ and an outer matrix function $F$ such that $\Psi=\Theta F$. 
Moreover, if $\Psi=\Theta^\circ F^\circ$, where $\Theta^\circ$ is inner and
$F^\circ$ is outer, then there exists a unitary matrix $\frak U$ such that
$\Theta^\circ=\Theta\frak U$ and $F^\circ=\frak U^*F$.

The above results can be found in the books [SNF1] and [Ni].

{\bf Badly approximable matrix functions and thematic factorizations.} 
Let $n\ge2$. An $n\times n$
matrix function $V$ is called {\it thematic} if it is unitary-valued and has the form
$$
\left(\begin{array}{cc}\bs{v}&\ov{\Theta}\end{array}\right),
$$
where $\bs{v}$ is an $n\times1$ inner and co-outer column function and $\Theta$ 
is an $n\times(n-1)$ inner and co-outer matrix function. It is natural to say that
a scalar function is thematic if it is constant and has modulus 1.

Let $\Phi$ be a matrix function in $L^\be(\mm_{m,n})$ such that 
$\|H_\Phi\|_{\text e}<\|H_\Phi\|$. If $F$ a best approximation of $\Phi$, then 
$\Phi-F$
admits a factorization
\beq
\label{2.1}
\Phi-F=W^*\left(\begin{array}{cc}t_0u&0\\0&\Psi\end{array}\right)V^*,
\end{equation}
where $V$ and $W^{\text t}$ are thematic matrix functions, $u$ is a unimodular 
scalar function (i.e., $|u(\z)|=1$ a.e. on $\T$) such that $T_u$ is Fredholm 
and $\ind T_u>0$, and 
$\Psi$ is a matrix function in $L^\be(\mm_{m-1,n-1})$ such that 
$\|\Psi\|_{L^\be}\le t_0$. This result was obtained in [PT]. Earlier the same fact 
was proved in [PY1] in the case \linebreak$\Phi\in(H^\be+C)(\mm_{m,n})$. Moreover, 
it was shown in [PT] that $\|H_\Psi\|_{\text e}\le\|H_\Phi\|_{\text e}$ and it was 
proved in [PY1] that the problem of
finding a superoptimal approximation of $\Phi$ reduces to the problem of finding
a superoptimal approximation of $\Psi$.

Clearly, the left-hand side of \rf{2.1} is a badly approximable function. 
Conversely, it follows from the results of [PY1] that
if $\Phi\in L^\be(\mm_{m,n})$ and $\Phi$ admits a factorization in
the form
$$
\Phi=W^*\left(\begin{array}{cc}su&0\\0&\Psi\end{array}\right)V^*,
$$
where $V$ and $W^{\text t}$ are thematic matrix functions, $s>0$, $u$ is a 
unimodular function such that $T_u$ is Fredholm and $\ind T_u>0$, and 
$\|\Psi\|_{L^\be}\le s$, then $\Phi$ is badly approximable and 
$s=t_0(\Phi)=\|H_\Phi\|$.

Suppose now that $\|H_\Phi\|_{\text e}<t_1$.
The inequality $\|H_\Psi\|_{\text e}\le\|H_\Phi\|_{\text e}$ proved in [PT] allows 
one to continue the diagonalization process and prove that if $l\le\min\{m,n\}$,
$\|H_\Phi\|_{\text e}<t_{l-1}$ and $t_{l-1}>t_l$, then for any matrix function
$F\in\bs{\O}_{l-1}$ (the sets $\bs{\O}_j$ are defined in \S 1) the matrix 
function $\Phi-F$ admits a factorization
\beq
\label{2.2}
\Phi-F=
W_0^*\cdots W^*_{l-1} 
\left( \begin{array}{ccccccc}
t_0 u_0 & 0 & \cdots & 0 & 0 \\
0 & t_1 u_1 & \cdots & 0 & 0 \\
\vdots & \vdots & \ddots & \vdots & \vdots\\
0 & 0 & \cdots & t_{l-1}u_{l-1} & 0\\
0&0&\cdots&0&\Psi\end{array}
\right)
V^*_{l-1}\cdots V^*_0,
\end{equation}
where
$$
W_j=\left(\begin{array}{cc}I_j & 0 \\ 0 & \breve{W}_j \end{array} \right),\quad
V_j = \left( \begin{array}{cc} I_j & 0 \\ 0 & \breve{V}_j \end{array} \right),
\quad 1 \leq j \leq l-1,
$$
the $W^{\text t}_0, \breve{W}^{\text t}_j, V_0, \breve{V}_j$ are thematic 
matrix functions, the $u_j$ are unimodular functions such that $T_{u_j}$ is
Fredholm and $\ind T_{u_j}>0$, $\|\Psi\|_{L^\be}\le t_{l-1}$, and
$\|H_\Psi\|<t_{l-1}$. Factorizations of
the form \rf{2.2} are called {\it partial thematic factorizations}.

Finally, if $\|H_\Phi\|_{\text e}$ is less than the smallest nonzero superoptimal
singular value of $\Phi$ and $F$ is the unique superoptimal approximation of $\Phi$,
then $\Phi-F$ admits a factorization
\beq
\label{2.3}
\!\!\!\!\!\!\Phi-F=
W_0^*\cdots W^*_{\iota-1} 
\left( \begin{array}{ccccccc}
t_0 u_0 & 0 & \cdots & 0 & 0 \\
0 & t_1 u_1 & \cdots & 0 & 0 \\
\vdots & \vdots & \ddots & \vdots & \vdots\\
0 & 0 & \cdots & t_{\iota-1}u_{\iota-1} & 0\\
0&0&\cdots&0&0\end{array}
\right)
V^*_{\iota-1}\cdots V^*_0,
\end{equation}
where the $u_j$, $V_j$, $W_j$ are as above and $t_{\iota-1}$ is the smallest nonzero
superoptimal singular value of $\Phi$. Factorizations of the form \rf{2.3}
are called {\it thematic factorizations}. This result was obtained in [PT]. Earlier
the same fact was proved in [PY1] in the case $\Phi\in(H^\be+C)(\mm_{m,n})$.
Note that the lower right entry of the diagonal matrix function on the right-hand 
side of \rf{2.3} has size $(m-\iota)\times(n-\iota)$ and the numbers $m-\iota$
or $n-\iota$ may be zero.

Clearly, the left-hand side of \rf{2.3} is a very badly approximable function.
Conversely, it follows from the results of [PY1] that if a matrix function
admits a thematic factorization, it is very badly approximable.

A disadvantage of thematic factorizations is that a thematic factorization
may essentially depend on the choice of matrix functions $V_j$ and $W_j$ and
they are not uniquely determined by the given matrix function. Moreover,
with any nonzero superoptimal singular value $t_j(\Phi)$ we can associate
the {\it factorization index} \linebreak$k_j\df\ind T_{u_j}$ of a thematic (or 
partial thematic) factorization. It was shown in [PY1] that in the case of multiple 
superoptimal singular values even the indices $k_j$ are not uniquely determined by 
$\Phi$. We refer the reader to [PY2], [PT], and [AP1] for further results on the
thematic factorization indices $k_j$. In particular, in [AP1] it was proved that
it is always possible to choose a so-called monotone (partial) thematic factorization
and the indices of a monotone (partial) thematic factorization are uniquely
determined by $\Phi$. However, the matrix function $V_j$ and $W_j$ are still
essentially dependent of our choice.

That is why we introduce in this paper (partial) canonical factorizations
and we prove that they are ``essentially'' unique.

\

\setcounter{equation}{0}
\setcounter{section}{2}
\section{\bf Balanced unitary-valued matrix functions}

\

In [PY1] and [PT] thematic matrix functions and thematic factorizations played a 
crucial role to study superoptimal approximation. We consider here a more general
class of balanced unitary-valued matrix functions.

{\bf Definition.} Let $n$ be a positive integer and let $r$ be an integer
such that $0<r<n$. Suppose that $\U$ is an $n\times r$ inner and co-outer matrix 
function and $\Theta$ is an $n\times(n-r)$ inner and co-outer matrix function. If
the matrix function
$$
\V=\left(\begin{array}{cc}\U&\ov{\Theta}\end{array}\right)
$$
is unitary-valued, it is called an {\it$r$-balanced matrix function}. 
If $r=0$ or $r=n$, it is natural 
to say that an $r$-balanced matrix is a constant unitary matrix.
An $n\times n$ matrix function $\V$ is called {\it balanced} if it is $r$-balanced 
for some $r$, $0\le r\le n$.

Recall that {\it thematic} matrix functions are just 1-balanced according to this 
definition.

It is known (see [V]) that if $0<r<n$ and $\U$ is an $n\times r$ inner and co-outer
matrix function, then it has a balanced unitary completion, i.e., there exists an
$n\times(n-r)$ matrix function $\Theta$ such that 
$\left(\begin{array}{cc}\U&\ov{\Theta}\end{array}\right)$ is balanced. 
Moreover, such a completion is unique modulo a right constant unitary factor. 
We are going to study some interesting properties of balanced matrix functions and
we need a construction of the complementary matrix function $\Theta$. That is why 
we give the construction here.

We will see that balanced matrix functions have many nice properties which can
justify the term ``balanced''.

Note that the case $r=1$ was studied in [PY1]. However, it turns out that studying
the more general case of an arbitrary $r$ simplifies the approach given in
[PY1].

\begin{thm}
\label{t3.1}
Let $n$ be a positive integer, $0<r<n$, and let $\U$ be and $n\times r$ inner
and co-outer matrix function. Then the subspace {\em$\cL\df\Ker T_{\U^{\text t}}$}
has the form
$$
\cL=\Theta H^2(\C^{n-r}),
$$
where $\Theta$ is an inner and co-outer $n\times(n-r)$ matrix function such that
$\left(\begin{array}{cc}\U&\ov{\Theta}\end{array}\right)$ is balanced.
\end{thm}

\Pf Clearly, the subspace $\cL$ of $H^2(\C^n)$ is invariant under multiplication
by $z$. By the Beurling--Lax theorem (see \S 2), $\cL=\Theta H^2(\C^l)$
for some $l\le n$ and an $n\times l$ inner matrix function $\Theta$.

Let us first prove that $\Theta$ is co-outer. Suppose that 
$\Theta^{\text t}=\cal O F$ where $\cal O$ is an inner matrix function and
$F$ is an outer matrix function. Since $\Theta$ is inner, it is easy to see
that $\cal O$ has size $l\times l$ while $F$ has size $l\times n$. It follows that
${\cal O}^*\Theta^{\text t}=F$, and so $F^{\text t}=\Theta\ov{\cal O}$. Since both
$\Theta$ and $\ov{\cal O}$ take isometric values almost everywhere on $\T$, the 
matrix function $F^{\text t}$ is inner.

Let us show that 
\beq
\label{3.1}
F^{\text t}H^2(\C^l)\subset\cL.
\end{equation} 
First of all, it is easy to see
that $\U^{\text t}\Theta$ is the zero matrix function.
Let now \linebreak$f\in H^2(\C^l)$. We have
$$
\U^{\text t}F^{\text t}f=\U^{\text t}\Theta\ov{\cal O}f=0,
$$
which proves \rf{3.1}.

It follows from \rf{3.1} that 
$$
F^{\text t}H^2(\C^l)=\Theta\ov{\cal O}H^2(\C^l)\subset\Theta H^2(\C^l).
$$
Multiplying the last inclusion by $\Theta^*$, we have
$$
\ov{\cal O}H^2(\C^l)\subset H^2(\C^l),
$$
which implies that $\cal O$ is a constant unitary matrix function.

Let us now prove that $l=n-r$. First of all it is evident that the columns of
$\U(\z)$ are orthogonal to the columns of $\ov{\Theta}(\z)$ almost everywhere on
$\T$. Hence, the matrix function 
$\left(\begin{array}{cc}\U&\ov{\Theta}\end{array}\right)$ takes isometric values 
almost everywhere, and so $l\le n-r$.

To show that $l\ge n-r$, consider the functions $P_\cL C$, where $P_\cL$
is the orthogonal projection onto $\cL$ and $C$ is a constant function which we
identify with a vector in $\C^n$. Note that $C\perp\cL$ if and only if
the vectors $f(0)$ and $C$ are orthogonal in $\C^n$ for any $f\in\cL$. Let us prove
that
\beq
\label{3.2}
\dim\{f(0):~f\in\cL\}\ge n-r.
\end{equation}
Since $\U$ is co-outer, it is easy to see that $\rank\U(0)=r$. Without loss
of generality we may assume that 
$\U=\left(\begin{array}{c}\U_1\\\U_2\end{array}\right)$, where
$\U_1$ has size $(n-r)\times r$, $\U_2$ has size $r\times r$, and 
the matrix $\U_2(0)$ is invertible. Let now $K$ be an arbitrary 
vector in $\C^{n-r}$. Put
$$
f=(\det\U_2(0))^{-1}\det\U_2\left(\begin{array}{c}K\\
-(\U_2^{\text t})^{-1}\U_1^{\text t}K
\end{array}\right).
$$
Clearly, $f\in H^2(\C^n)$. We have
$$
\U^{\text t}f=\left(\begin{array}{cc}\U_1^{\text t}&\U_2^{\text t}\end{array}\right)f
=(\det\U_2(0))^{-1}\det\U_2(\U_1^{\text t}K-
\U_2^{\text t}(\U_2^{\text t})^{-1}\U_1^{\text t}K)=0,
$$
and so $f\in\cL$. On the other hand, it is easy to see
that
$$
f(0)=\left(\begin{array}{c}K\\-(\U_2^{\text t}(0))^{-1}\U_1^{\text t}(0)K
\end{array}\right)
$$
and since $K$ is an arbitrary vector in $\C^{n-r}$, this proves \rf{3.2}.

We have already observed that
$$
\{f(0):~f\in\cL\}=\C^n\ominus\{C\in\C^n:~P_\cL C=0\},
$$
and so it follows from \rf{3.2} that
$$
\dim\{P_\cL C:~C\in\C^n\}\ge n-r.
$$
It is easy to see that for $C\in\C^n$ we have
$$
P_\cL C=\Theta\pp_+\Theta^*C=\Theta(\Theta^*(0))C.
$$
Clearly, $\Theta(\Theta^*(0))C$ belongs to the linear span of the columns of 
$\Theta$. This completes the proof of the fact that $l=n-r$ and proves that 
$\left(\begin{array}{cc}\U&\ov{\Theta}\end{array}\right)$ is a balanced
matrix function. $\bl$

It is well known (see [V]) that given an inner and co-outer matrix function $\U$,
the balanced completion $\Theta$ is unique modulo a right constant unitary factor.
Indeed, if $\left(\begin{array}{cc}\U&\ov{\Theta}_1\end{array}\right)$ is
another balanced matrix function, then clearly 
$$
\Theta_1H^2(\C^{n-r})\subset\cL=\Theta H^2(\C^{n-r}).
$$
It follows (see [Ni]) that $\Theta_1=\Theta{\cal O}$ for an inner matrix function 
$\cal O$. Clearly, $\cal O$ has size $(n-r)\times(n-r)$. Hence, 
$\Theta_1^{\text t}={\cal O}^{\text t}\Theta^{\text t}$, ${\cal O}^{\text t}$ is
inner, and since $\Theta_1$ is co-outer, it follows that $\cal O$ is a
unitary constant.

Next, we are going to study the property of analyticity of minors of balanced 
matrix functions. Let
$\V=\left(\begin{array}{cc}\U&\ov{\Theta}\end{array}\right)$ be an $r$-balanced
$n\times n$ matrix function. We are going to study its minors 
$\V_{\i_1\cdots\i_k,\j_1\cdots\j_k}$ of order $k$, i.e., the determinants of
the submatrix of $\V$ with rows $\i_1,\cdots,\i_k$ and columns 
$\j_1,\cdots,\j_k$. Here \linebreak$1\le\i_1<\cdots<\i_k\le n$ 
and $1\le\j_1<\cdots<\j_k\le n$.
By a {\it minor of $\V$ on the first $r$ columns} we mean a minor
$\V_{\i_1\cdots\i_k,\j_1\cdots\j_k}$ with $k\ge r$ and $\j_1=1,\,\cdots,\,\j_r=r$.
Similarly, by a {\it minor of $\V$ on the last $n-r$ columns} we mean a minor
$\V_{\i_1\cdots\i_k,\j_1\cdots\j_k}$ with $k\ge n-r$ and 
$\j_{k-n+r+1}=r+1,\,\cdots,\,\j_k=n$.

It was proved in [PY1] (Theorem 1.1) that for a thematic (i.e., 1-balanced) 
matrix function all minors on the first column are in $H^\be$. The following
result generalizes that theorem and makes it more symmetric.

\begin{thm}
\label{t3.2}
Let $\V$ be an $r$-balanced matrix function of size $n\times n$. Then all 
minors of $\V$ on the first $r$ columns are in $H^\be$ while all minors of $\V$
on the last $n-r$ columns are in $\ov{H^\be}$.
\end{thm}

\Pf It is easy to see that it is sufficient to prove the theorem for the minors
on the first $r$ columns of $\V$. To deduce the second assertion of the theorem, 
we can consider the matrix function $\ov{\V}$ and rearrange its columns to make it
\linebreak$(n-r)$-balanced.

Denote by $\U_1,\cdots,\U_r$ and $\Theta_1,\cdots,\Theta_{n-r}$ the columns of $\U$
and $\Theta$. In the proof of Theorem \ref{t3.1} we have observed that for any
constant $C\in\C^n$ we have
$$
P_\cL C=\Theta\Theta^*(0)C
$$
and
$$
\dim\{P_\cL C:~C\in\C^n\}=\dim\{\Theta^*(0)C:~C\in\C^n\}=n-r.
$$
It follows that there exist $C_1,\cdots,C_{n-r}\in\C^n$ such that
\beq
\label{3.3}
\Theta_j=P_\cL C_j,\quad1\le j\le n-r.
\end{equation}
If $1\le d\le n-r$ and $1\le j_1<j_2<\cdots<j_d\le n-r$, we consider the 
vector function
$$
\U_1\wedge\cdots\wedge\U_r\wedge\ov{\Theta}_{j_1}\wedge\cdots\wedge\ov{\Theta}_{j_d}
$$
whose $\left(\begin{array}{c}n\\r+d\end{array}\right)$ components are the 
minors of order $r+d$ of the matrix function
$\left(\begin{array}{cccccc}
\U_1&\cdots&\U_r&\ov{\Theta}_{j_1}&\cdots&\ov{\Theta}_{j_d}\end{array}\right)$.

It follows from \rf{3.3} that
$$
\Theta_j=C_j-P_{\cL^\perp}C_j,
$$
where $P_{\cL^\perp}$ is the orthogonal projection onto 
$\cL^\perp=\clos\Range T_{\ov{\U}}$. We have
$$
\U_1\wedge\cdots\wedge\U_r\wedge\ov{\Theta}_{j_1}\wedge\cdots\wedge\ov{\Theta}_{j_d}=
\U_1\wedge\cdots\wedge\U_r\wedge
(\ov{C}_{j_1}-\ov{P_{\cL^\perp}C_{j_1}})\wedge\cdots
\wedge(\ov{C}_{j_d}-\ov{P_{\cL^\perp}C_{j_d}}).
$$
The components of this vector belong to $L^\be$ and can be approximated in $L^{2/d}$
by vector functions of the form
\beq
\label{3.4}
\U_1\wedge\cdots\wedge\U_r\wedge
(\ov{C}_{j_1}-\ov{g}_{j_1})\wedge\cdots
\wedge(\ov{C}_{j_d}-\ov{g}_{j_d}),
\end{equation}
where $g_{j_1},\cdots,g_{j_d}\in\Range T_{\ov{\U}}$. Hence, it is sufficient to
prove that the components of \rf{3.4} belong to $H^{2/d}$. Let 
$g_{j_l}=\pp_+\ov{\U}f_l$ for $f_l\in H^2(\C^r)$, $1\le l\le d$. We have
\begin{eqnarray*}
&&\U_1\wedge\cdots\wedge\U_r\wedge
(\ov{C}_{j_1}-\ov{g}_{j_1})\wedge\cdots
\wedge(\ov{C}_{j_d}-\ov{g}_{j_d})\\
&=&\U_1\wedge\cdots\wedge\U_r\wedge
(\ov{C}_{j_1}-\ov{\pp_+\ov{\U}f_1})\wedge\cdots
\wedge(\ov{C}_{j_d}-\ov{\pp_+\ov{\U}f_d})\\
&=&\U_1\wedge\cdots\wedge\U_r\wedge
(\ov{C}_{j_1}-\U\bar f_1+\ov{\pp_-\ov{\U}f_1})\wedge\cdots
\wedge(\ov{C}_{j_d}-\U\bar f_d+\ov{\pp_-\ov{\U}f_d}).
\end{eqnarray*}
Clearly, almost everywhere on $\T$ the vectors $\U(\z)\ov{f_l(\z)}$ are linear 
combinations of $\U_1(\z),\cdots,\U_r(\z)$. Therefore if we expand the above wedge
product using the multilinearity of $\wedge$, all terms containing $\U\bar f_l$,
give zero contribution. Thus we have
\begin{eqnarray*}
&&\U_1\wedge\cdots\wedge\U_r\wedge
(\ov{C}_{j_1}-\U\bar f_1+\ov{\pp_-\ov{\U}f_1})\wedge\cdots
\wedge(\ov{C}_{j_d}-\U\bar f_d+\ov{\pp_-\ov{\U}f_d})\\
&=&\U_1\wedge\cdots\wedge\U_r\wedge
(\ov{C}_{j_1}+\ov{\pp_-\ov{\U}f_1})\wedge\cdots
\wedge(\ov{C}_{j_d}+\ov{\pp_-\ov{\U}f_d})\in H^{2/d}.\quad\bl
\end{eqnarray*}

The following immediate consequence of Theorem \ref{t3.2} was obtained in [PY1]
(Theorem 1.2) for $r=1$ by another method.

\begin{cor}
\label{t3.3}
Let $\V$ be a balanced matrix function. Then $\det\V$ is a constant function 
of modulus 1.
\end{cor}

We shall need another nice property of balanced matrix functions that was 
obtained in [Pe4]. Namely, it was proved in Lemma 6.2 of [Pe4] that for a balanced 
matrix function $\V$ of size $n\times n$ the Toeplitz operator 
$T_\V:H^2(\C^n)\to H^2(\C^n)$ has trivial kernel and dense range. If we apply that
result to $\ov\V$ with rearranged columns, we see that the Toeplitz operator
$T_{\ov{\V}}$ also has trivial kernel and dense range.

Finally, we obtain in this section an analog of Lemma 1.5 of [PY1] where the case
of thematic matrix functions was considered.

\begin{thm}
\label{t3.4}
Let $0<r\le\min\{m,n\}$ and let $\V$ and {\em$\W^{\text t}$} be $r$-balanced
matrix functions of sizes $n\times n$ and $m\times m$ respectively. Then
$$
\W H^\be(\mm_{m,n})\V\bigcap
\left(\begin{array}{cc}0&0\\0&L^\be(\mm_{m-r,n-r})\end{array}\right)=
\left(\begin{array}{cc}0&0\\0&H^\be(\mm_{m-r,n-r})\end{array}\right).
$$
\end{thm}

The proof of Theorem \ref{t3.4} is exactly the same as the proof of its special
case Lemma 1.5 of [PY1].

\

\setcounter{equation}{0}
\setcounter{section}{3}
\section{\bf Best approximation and partial canonical factorizations}

\

For a matrix function $\Phi$ in $L^\be(\mm_{m,n})$ satisfying the condition
$\|H_\Phi\|_{\text e}<\|H_\Phi\|$ and a best approximation
$F$ of $\Phi$ by bounded analytic matrix functions we obtain a so-called
partial canonical factorization of $\Phi-F$. We characterize badly approximable
functions satisfying this condition in terms of partial canonical
factorizations. To this end we begin this section with the study of
the minimal invariant subspace of multiplication by $z$ on $H^2(\C^n)$
that contains all maximizing vectors of $H_\Phi$.

\begin{thm}
\label{t4.1}
Let $\Phi\in L^\be(\mm_{m,n})$ and 
{\em$\|H_\Phi\|_{\text e}<\|H_\Phi\|$}. Let $\M$ be the minimal invariant
subspace of multiplication by $z$ on $H^2(\C^n)$ that contains all maximizing 
vectors of $H_\Phi$. Then
\beq
\label{4.1}
\M=\U H^2(\C^r),
\end{equation}
where $r$ is the number of superoptimal singular values of $\Phi$ equal to 
$\|H_\Phi\|$, $\U$ is an inner and co-outer $n\times r$ matrix function.
\end{thm}

\Pf Consider first the case $m=n$. Without loss of generality we may assume that
$\|H_\Phi\|=1$. It follows from the results of \S 3 of  [AP2]
that there exists a unitary interpolant 
$\cal U$ of $\Phi$ (i.e., a unitary-valued matrix function $\cal U$  satisfying 
$\hat{\cal U}(j)=\hat\Phi(j)$, $j<0$) such that the Toeplitz operator $T_{\cal U}$
is Fredholm and each such unitary interpolant has precisely $r$ negative
Wiener--Hopf indices. Consider a Wiener--Hopf factorization of $\cal U$
\beq
\label{4.2}
{\cal U}=Q^*_2\left(\begin{array}{cccc}z^{d_1}&0&\cdots&0\\
0&z^{d_2}&\cdots&0\\\vdots&\vdots&\ddots&\vdots\\0&0&\cdots&z^{d_n}
\end{array}\right)Q_1^{-1},
\end{equation}
where $Q_1$ and $Q_2$ are matrix functions invertible in $H^2(\mm_{n,n})$, and
$d_1\le d_2\cdots\le d_n$. Since $\cal U$ has $r$ negative Wiener--Hopf indices,
we have 
$$
d_1\le\cdots\le d_r<0\le d_{r+1}\le\cdots\le d_n.
$$

Clearly, $H_\Phi=H_{\cal U}$. It is also easy to see that a nonzero function
$f\in H^2(\C^n)$ is a maximizing vector of $H_\Phi$ if an only if 
$f\in\Ker T_{\cal U}$. It is well known and it is easy to see from (\ref{4.2}) that
\beq
\label{4.3}
\Ker T_{\cal U}=
\left\{Q_1\left(\begin{array}{c}q_1\\\vdots\\q_r\\0\\\vdots\\0\end{array}
\right):~q_j\in\PP_+,~\deg q_j<-d_j,~1\le j\le r
\right\}.
\end{equation}
Here we denote by $\PP_+$ the set of analytic polynomials.
Since $\M$ is the minimal invariant subspace of multiplication by $z$ that
contains $\Ker T_{\cal U}$, it follows from (\ref{4.3}) that
\beq
\label{4.4}
\M=\clos_{H^2(\C^n)}
\left\{Q_1\left(\begin{array}{c}q_1\\\vdots\\q_r\\0\\\vdots\\0\end{array}
\right):~q_j\in\PP_+,~1\le j\le r
\right\}.
\end{equation}
Since $Q_1(\z)$ is an invertible matrix for all $\z\in\dd$, 
it follows easily from \rf{4.4} that $\dim\{f(\z):~f\in\M\}=r$ for all
$\z\in\dd$. Therefore the $z$-invariant subspace $\M$ has the form
$\M=\U H^2(\C^r)$, where $\U$ is an $n\times r$ inner matrix function (see [Ni]).
It remains to prove that $\U$ is co-outer.

Denote by $Q_\heartsuit$ the matrix function obtained from $Q_1$ by deleting 
the last $n-r$ columns. It is easy to see that
$\U$ is an inner part of $Q_\heartsuit$. Let $Q_\heartsuit=\U F$, where $F$ is an 
$r\times r$ outer matrix function. 

Denote by $Q_\spadesuit$ the matrix function obtained from $Q_1^{-1}$ by deleting 
the last $n-r$ rows. Clearly, $Q_\spadesuit(\z)Q_\heartsuit(\z)=I_r$
for almost all $\z\in\T$.

We have 
$$
I_r=Q_\spadesuit Q_\heartsuit=Q_\spadesuit\U F,
$$
and so
$$
I_r=F^{\text t}\U^{\text t}Q_\spadesuit^{\text t}.
$$
Both $F^{\text t}$ and $\U^{\text t}Q_\spadesuit^{\text t}$ are $r\times r$
matrix functions. Hence,
$$
I_r=\U^{\text t}Q_\spadesuit^{\text t}F^{\text t}.
$$
It follows that $\U^{\text t}$ is outer, and so $\U$ is co-outer.

Consider now the case $m<n$. Let $\Phi_\#$ be the matrix function obtained 
from $\Phi$ by adding $n-m$ zero rows. It is easy to see that the Hankel operators
$H_\Phi$ and $H_{\Phi_\#}$ have the same maximizing vectors. This reduces
the problem to the case $m=n$.

Finally, assume that $m>n$. Let $\Phi_\flat$ be the matrix function obtained 
from $\Phi$ by adding $m-n$ zero columns. It is easy to see that $f$ is a maximizing 
vector of $H_{\Phi_\flat}$  if and only if it can be obtained from a maximizing
vector of $H_\Phi$ by adding $m-n$ zero coordinates. Let $\M_\flat$ be the 
minimal invariant subspace of multiplication by $z$ on $H^2(\C^m)$ that contains
all maximizing vectors of $H_{\Phi_\flat}$. Clearly, the number of superoptimal 
singular values of $\Phi_\flat$ equal to 1 is still $r$. Therefore there exists
an $m\times r$ inner and co-outer matrix function $\U_\flat$ such that
$$
\M_\flat=\U_\flat H^2(\C^r).
$$
It is easy to see that the last $m-n$ rows of $\U_\flat$ are zero. Denote by $\U$
the matrix function obtained from $\U_\flat$ by deleting the last $m-n$ zero
rows. Obviously, $\U$ is an inner and co-outer $n\times r$ matrix function
and $\M=\U H^2(\C^r)$. $\bl$

We need the following result.

\begin{lem}
\label{t4.2}
Suppose that $\Phi$ satisfies the hypotheses of Theorem \ref{t4.1} and $\M$ is given 
by {\em\rf{4.1}}. If $\|\Phi\|_{L^\be}=\|H_\Phi\|$ and $f$ is a nonzero vector 
function in $\M$, then $f(\z)$ is a maximizing vector of $\Phi(\z)$ for almost all 
$\z\in\T$.
\end{lem}

\Pf As we have noted in \S 2, if $f$ is a maximizing vector of $H_\Phi$, then
$f(\z)$ is a maximizing vector of $\Phi$ for almost all $\z\in\T$ and 
$\|\Phi(\z)\|_{\mm_{m,n}}=\|H_\Phi\|$ almost everywhere. Without loss of generality
we may assume that $\|H_\Phi\|=1$.

Let $L$ be the set of vector functions of the form 
$$
q_1g_1+\cdots+q_Mg_M,
$$ 
where $q_j\in\PP_+$ and the $g_j$ are maximizing vectors of $H_\Phi$. 
By definition, $\M$ is the norm closure of $L$. Since
the $g_j(\z)$ are maximizing vectors of $\Phi(\z)$ for almost all $\z\in\T$ 
(see \S 2), it follows that for $g\in L$, $g(\z)$ is a maximizing vector of 
$\Phi(\z)$ almost everywhere on $\T$. Let $\{f_j\}$ be a sequence of vector functions
in $L$ that converges to $f\in\M$ in $H^2(\C^n)$. Clearly,
\begin{eqnarray*}
\int_\T\|\Phi(\z)f(\z)\|_{\C^m}^2d\m(\z)&=&
\lim_{j\to\be}\int_\T\|\Phi(\z)f_j(\z)\|_{\C^m}^2d\m(\z)\\
&=&\lim_{j\to\be}\int_\T\|f_j(\z)\|_{\C^m}^2d\m(\z)\\
&=&\int_\T\|f(\z)\|_{\C^m}^2d\m(\z),
\end{eqnarray*}
and since obviously, $\|\Phi(\z)f(\z)\|_{\C^m}\le\|f(\z)\|_{\C^m}$ almost
everywhere on $\T$, it follows that $f(\z)$ is a maximizing vector of $\Phi(\z)$
for almost all $\z\in\T$. $\bl$

Again, suppose that $\Phi$ is as in Theorem \ref{t4.1}. Obviously, a vector function 
$f$ in $H^2(\C^n)$ is a maximizing vector of $H_\Phi$ if and only if 
$g\df H_\Phi f\in H^2_-(\C^m)$ 
is a maximizing vector of $H_\Phi^*$. Let us show that a vector function $g$ in 
$H^2_-(\C^m)$ is a maximizing vector of $H_\Phi^*$ if and only if 
$\bar z\bar g\in H^2(\C^m)$
is a maximizing vector of $H_{\Phi^{\text t}}$. Indeed, assume that $\|H_\Phi\|=1$.
Then $g$ is a maximizing vector of $H_\Phi^*$ if and only if
$$
\|H_\Phi^*g\|_2=\|\pp_+\Phi^*g\|_2=\|g\|_2.
$$
Clearly, this is equivalent to the equality
$$
\|\pp_-\Phi^{\text t}\bar z\bar g\|_2=\|\bar z\bar g\|_2,
$$
which means that $\bar z\bar g$ is a maximizing vector of $H_{\Phi^{\text t}}$.

It is easy to see that the matrix functions $\Phi$ and $\Phi^{\text t}$ have
the same superoptimal singular values. Let $\cal N$ be the minimal invariant
subspace of multiplication by $z$ on $H^2(\C^m)$ that contains all maximizing vectors
of $H_{\Phi^{\text t}}$. By Theorem \ref{t4.1}, there exists an inner and co-outer 
matrix function $\O\in H^\be(\mm_{m,r})$ such that
$$
\cal N=\O H^2(\C^r).
$$
By Theorem \ref{t3.1}, $\U$ and $\O$ have balanced completions, i.e., 
there exist inner and co-outer matrix functions 
$\Theta\in H^\be(\mm_{n,n-r})$ and $\Xi\in H^\be(\mm_{m,m-r})$ such that
\beq
\label{4.5}
\V\df\left(\begin{array}{cc}\U&\ov{\Theta}\end{array}\right)\quad\mbox{and}\quad
\W^{\text t}\df\left(\begin{array}{cc}\O&\ov{\Xi}\end{array}\right)
\end{equation}
are unitary-valued matrix functions.

\begin{thm}
\label{t4.3}
Let $\Phi\in L^\be(\mm_{m,n})$ and {\em$\|H_\Phi\|_{\text e}<t_0=\|H_\Phi\|$}. 
Let $r$ be the number of superoptimal singular values of $\Phi$ equal to 
$t_0$. Suppose that $F$ is a best approximation
of $\Phi$ by analytic matrix functions. Then $\Phi-F$ admits a factorization
of the form
\beq
\label{4.6}
\Phi-F=\W^*\left(\begin{array}{cc}t_0U&0\\0&\Psi
\end{array}\right)\V^*,
\end{equation}
where $\V$ and $\W$ are given by {\em\rf{4.5}}, $U$ is an $r\times r$ unitary-valued 
very badly approximable matrix function such that {\em$\|H_U\|_{\text e}<1$},
and $\Psi$ is a matrix-function in
$L^\be(\mm_{m-r,n-r})$ such that $\|\Psi\|_{L^\be}\le t_0$ and
$\|H_\Psi\|=t_r(\Phi)<\|H_\Phi\|$. Moreover, $U$ is uniquely determined by the 
choice of $\U$ and $\O$, and does not depend on the choice of $F$.
\end{thm}

\Pf Without loss of generality we may assume that $\|H_\Phi\|=1$. It follows from 
Lemma \ref{t4.2} that the columns of $\U(\z)$ are maximizing vectors of 
$\Phi(\z)-F(\z)$ for almost all $\z\in\T$. Similarly, the columns of $\O(\z)$ are 
maximizing vectors of $\Phi^{\text t}(\z)-F^{\text t}(\z)$ almost everywhere on $\T$.

We need two elementary lemmas.

\begin{lem}
\label{t4.4}
Let $A\in\mm_{m,n}$ and $\|A\|=1$. Suppose that $v_1,\cdots,v_r$ is an 
orthonormal family of maximizing vectors of $A$ and $w_1,\cdots,w_r$ is an 
orthonormal family of maximizing vectors of $A^{\text t}$. Then{\em
$$
\left(\begin{array}{ccc}w_1&\cdots&w_r\end{array}\right)^{\text t}A
\left(\begin{array}{ccc}v_1&\cdots&v_r\end{array}\right)
$$}
is a unitary matrix.
\end{lem}

\begin{lem}
\label{t4.5}
Let $A$ be a matrix in $\mm_{m,n}$ such that $\|A\|=1$ and $A$ has the form
$$
A=\left(\begin{array}{cc}A_{11}&A_{12}\\A_{21}&A_{22}\end{array}\right),
$$
where $A_{11}$ is a unitary matrix. Then $A_{12}$ and $A_{21}$ are the zero matrices.
\end{lem}

Both lemmas are obvious. Let us complete the proof of Theorem \ref{t4.3}.

Consider the matrix function
$$
\left(\begin{array}{cc}U&X\\Y&\Psi\end{array}\right)\df\W(\Phi-F)\V.
$$
Here $U\in L^\be(\mm_{r,r})$, $X\in L^\be(\mm_{r,n-r})$, $Y\in L^\be(\mm_{m-r,r})$,
and $\Psi\in L^\be(\mm_{m-r,n-r})$. If we apply Lemma \ref{t4.4} to the matrices 
$(\Phi-F)(\z)$, $\z\in\T$, and the columns of $\U(\z)$ and $\O(\z)$, we see 
$U=\O^{\text t}(\Phi-F)\U$ is 
unitary-valued. By Lemma \ref{t4.5}, $X$ and $Y$ are the zero matrix functions which 
proves \rf{4.6}. 

Let us prove that $\|H_U\|_{\text e}<1$. Since 
\beq
\label{4.7}
U=\O^{\text t}(\Phi-F)\U,
\end{equation}
we have
\begin{eqnarray*}
\|H_U\|_{\text e}&=&\dist_{L^\be}\big(U,(H^\be+C)(\mm_{r,r})\big)\\
&=&\dist_{L^\be}\big(\O^{\text t}\Phi\U,(H^\be+C)(\mm_{r,r})\big)\\
&\le&\dist_{L^\be}\big(\Phi,(H^\be+C)(\mm_{m,n})\big)=\|H_\Phi\|_{\text e}<1.
\end{eqnarray*}

Now it is turn to show that $U$ is very badly approximable. 
Denote by $\cL$ the minimal invariant subspace of multiplication by $z$
that contains all maximizing vectors of $H_U$. Suppose that
$f$ is a maximizing vector of $H_\Phi$. Then $f=\U g$ for some $g\in H^2(\C^r)$.
Clearly, $H_\Phi f$ is a maximizing vector of $H_\Phi^*$.
As we have mentioned in \S 2, $H_\Phi f=(\Phi-F)f$. Hence, 
$\bar z\ov{\Phi-F}\bar f$
is a maximizing vector of $H_{\Phi^{\text t}}$. Therefore 
$\bar z\ov{\Phi-F}\bar f\in\O H^2(\C^r)$, and so 
$$
\O^{\text t}(\Phi-F)f=\O^{\text t}(\Phi-F)\U g=U g\in H^2_-(\C^r).
$$
It follows that $g$ is a maximizing vector of $H_U$ and $\|H_U\|=1$. Therefore
$$
\U H^2(\C^r)\subset\U\cL.
$$
Hence, $\cL=H^2(\C^r)$, and by Theorem \ref{t4.1}, $t_0(U)=\cdots=t_{r-1}(U)=1$.
It follows that $U$ is very badly approximable. Hence, the zero matrix
function is the only best approximation of $U$ by analytic matrix functions.

This uniqueness property together with \rf{4.7} implies that $U$ does not depend on 
the choice of the best approximation $F$.

It is evident from \rf{4.6} that $\|\Psi\|_{L^\be}\le1$. It remains to prove that
$\|H_\Psi\|=t_r(\Phi)$.

Suppose that $F_{\$}$ is another best approximation of $\Phi$ by bounded analytic 
matrix functions. Then as we have already proved, $\Phi-F_{\$}$ can be represented
as
$$
\Phi-F_{\$}=\W^*\left(\begin{array}{cc}U&0\\0&\Psi_{\$}
\end{array}\right)\V^*,
$$
where $\Psi_{\$}$ is a matrix function in $L^\be(\mm_{m-r,n-r})$ such that
$\|\Psi_{\$}\|_{L^\be}\le1$. Clearly, $s_j((\Phi-F_{\$})(\z))=1$,
$0\le j\le r-1$,  and $s_r((\Phi-F_{\$})(\z))=\|\Psi_{\$}(\z)\|$ for almost all 
$\z\in\T$.

By Theorem \ref{t3.4}, a matrix function $G\in H^\be(\mm_{m,n})$ is a best 
approximation of $\Phi$ if and only if there exists $Q\in H^\be(\mm_{m-r,n-r})$ such 
that $\|\Psi-Q\|_{L^\be}\le1$ and
$$
\Phi-G=\W^*\left(\begin{array}{cc}U&0\\0&\Psi-Q
\end{array}\right)\V^*.
$$
This proves that $\|H_\Psi\|=t_r(\Phi)$. $\bl$

{\bf Remark.} It can be shown easily that if $U$ is an $r\times r$ unitary-valued 
very badly approximable matrix function and $\|H_U\|_{\text e}<1$, then the Toeplitz 
operator $T_U$ is Fredholm.
Indeed, by Theorem 3.1 of [PY1], if $\Phi$ is a very badly approximable function
in $(H^\be+C)(\mm_{m,n})$ and $\rank\Phi(\z)=m$ almost everywhere on $\T$, then the
Toeplitz operator $T_{z\Phi}:H^2(\C^n)\to H^2(\C^m)$ has dense range. The results
of [PT] show that the proof given in [PY1] also works in the more general case 
when $\Phi\in L^\be(\mm_{m,n})$ and $\|H_\Phi\|_{\text e}<\|H_\Phi\|$. Hence,
$T_{zU}$ has dense range. Then the operator $H^*_{\bar zU^*}H_{\bar zU^*}$
is unitarily equivalent to the restriction of $H^*_{zU}H_{zU}$ to the orthogonal
complement to the subspace
$$
\{f\in H^2(\C^r):~\|H_{zU}f\|_2=\|f\|_2\}
$$
(see [Pe2]). Since $\|H_U\|_{\text e}<1$, this subspace is finite-dimensional, and
so 
$$
\|H_{U^*}\|_{\text e}=\|H_{\bar zU^*}\|_{\text e}=\lim_{j\to\be}s_j(H_{\bar zU^*})=
\lim_{j\to\be}s_j(H_{zU})=\|H_{zU}\|_{\text e}=\|H_U\|_{\text e}.
$$
The result follows now from the well-known fact that for a unitary-valued function 
$U$ the conditions $\|H_U\|_{\text e}<1$ and $\|H_{U^*}\|_{\text e}<1$ are equivalent
to the fact that $T_U$ is Fredholm (see e.g., \S 1 of [AP2] for some comments).

Factorizations of the form \rf{4.6} with $\Psi$ satisfying
$$
\|\Psi\|_{L^\be}\le t_0\quad\mbox{and}\quad\|H_\Psi\|<t_0
$$
form a special class of {\it partial canonical 
factorizations}. The matrix function $\Psi$ is called the {\it residual entry}
of the partial canonical factorization.
The notion of a partial canonical factorization in the general
case will be defined in \S 7.

The following theorem together with Theorem \ref{t4.3} gives a characterization
of the badly approximable matrix functions $\Phi$ satisfying the condition
$\|H_\Phi\|_{\text e}<\|H_\Phi\|$. 

\begin{thm}
\label{t4.6}
Let $\Phi\in L^\be(\mm_{m,n})$ and {\em$\|H_\Phi\|_{\text e}<\|H_\Phi\|$}.
Suppose that $\Phi$ admits a representation of the form
$$
\Phi=\W^*\left(\begin{array}{cc}\s U&0\\0&\Psi
\end{array}\right)\V^*,
$$
where $\s>0$, $\V$ and {\em$\W^{\text t}$} are $r$-balanced matrix functions, 
$U$ is a very badly approximable unitary-valued $r\times r$ matrix function such that
{\em$\|H_U\|_{\text e}<1$}, and $\|\Psi\|_{L^\be}\le\s$. Then
$\Phi$ is badly approximable and $t_0(\Phi)=\cdots=t_{r-1}(\Phi)=\s$.
\end{thm}

\Pf Suppose that $\V$ and $\W$ are given by \rf{4.5}.
Let $g\in H^2(\C^r)$ be a maximizing vector of $H_U$. Then it is easy to see that
$$
\|H_\Phi\U g\|_2=\s\|\U g\|,
$$
while
$$
\|H_\Phi\|\le\|\Phi\|_{L^\be}=\s.
$$
Hence, $\|H_\Phi\|=\s$, $\U g$ is a maximizing vector of $H_\Phi$, and $\Phi$ is 
badly approximable.

Since $U$ is very badly approximable, it follows from Theorem \ref{t4.1} that
the minimal invariant subspace of multiplication by $z$ on $H^2(\C^r)$ that
contains all maximizing vectors of $H_U$ is the space $H^2(\C^r)$ itself.
Let $\M$ be the minimal invariant subspace of multiplication by $z$ on $H^2(\C^n)$
that contains all maximizing vectors of $H_\Phi$. Since $\U$ is co-outer, it follows 
that the matrix $\U(\z)$ has rank $r$ for all $\z\in\dd$.
Hence,
$$
\dim\{f(\z):~f\in\M\}\ge r\quad \mbox{for all}\quad\z\in\dd.
$$
It follows now from Theorem \ref{t4.1} that $t_0(\Phi)=\cdots=t_{r-1}(\Phi)$. $\bl$

We also need the following version of the converse to Theorem \ref{t4.3}.

\begin{thm}
\label{t4.7}
Let $\Phi\in L^\be(\mm_{m,n})$ and {\em$\|H_\Phi\|_{\text e}<\|H_\Phi\|$}.
Suppose that $\Phi$ admits a representation of the form
$$
\Phi=\W^*\left(\begin{array}{cc}\s U&0\\0&\Psi
\end{array}\right)\V^*,
$$
where $\s>0$, $\V$ and {\em$\W^{\text t}$} are $r$-balanced matrix functions
of the form {\em\rf{4.5}}, 
$U$ is a very badly approximable unitary-valued $r\times r$ matrix function such that
{\em$\|H_U\|_{\text e}<1$}, and $\|H_\Psi\|<\s$. Then $\U H^2(\C^r)$ is the 
minimal invariant subspace of multiplication by $z$ on $H^2(\C^n)$ that contains all 
maximizing vectors of $H_\Phi$ and $\O H^2(\C^r)$ is the minimal
invariant subspace of multiplication by $z$ on $H^2(\C^m)$ that contains all 
maximizing vectors of {\em$H_{\Phi^{\text t}}$}.
\end{thm}

\Pf By Theorem \ref{t3.4}, we may assume without loss of generality that
\linebreak$\|\Psi\|_{L^\be}<s$. We need the following lemma.

\begin{lem}
\label{t4.8}
Suppose that $\Phi$ satisfies the hypotheses of Theorem \ref{t4.7}.
A function $f$ in $H^2(\C^n)$ is a maximizing vector of $H_\Phi$
if and only if $f=\U g$, where $g\in H^2(\C^r)$ and $g$ is a maximizing vector
of $H_U$.
\end{lem}

Let us first complete the proof of Theorem \ref{t4.7}. Since $U$ is very badly
approximable, we have $t_0(U)=\cdots=t_{r-1}(U)=1$. By Theorem \ref{t4.1},
the minimal invariant subspace of multiplication by $z$ on $H^2(\C^r)$ that contains
all maximizing vectors of $H_U$ is $H^2(\C^r)$. It follows now from
Lemma \ref{t4.8} that the 
minimal invariant subspace of multiplication by $z$ on $H^2(\C^n)$ that contains all 
maximizing vectors of $H_\Phi$ is $\U H^2(\C^r)$. To complete the proof,
we can apply this result to the matrix function $\Phi^{\text t}$. $\bl$

{\bf Proof of Lemma \ref{t4.8}.} First of all, by Theorem \ref{t4.6},
$\|H_\Phi\|=\s$. Without loss of generality we may assume that $\s=1$.
It has been proved in the proof of Theorem \ref{t4.6} that if $g$ is a maximizing
vector of $H_U$, then $\U g$ is a maximizing vector of $H_\Phi$.

Suppose now that $f$ is a maximizing vector of $H_\Phi$.
We have
\begin{eqnarray*}
\Phi f&=&\W^*\left(\begin{array}{cc}U&0\\0&\Psi\end{array}\right)\V^*f\\
&=&\W^*\left(\begin{array}{cc}U&0\\0&\Psi\end{array}\right)
\left(\begin{array}{c}\U^*f\\\Theta^{\text t}f\end{array}\right)\\
&=&\W^*\left(\begin{array}{c}U\U^*f\\\Psi\Theta^{\text t}f\end{array}\right).
\end{eqnarray*}

Since $\W^*$ is unitary-valued and $\|\Psi\|_{L^\be}<1$, it follows that
$\Theta^{\text t}f=0$. Put $g=\U^*f\in L^2(\C^r)$. We have
$$
f=\V\V^*f=\V\left(\begin{array}{c}\U^*\\\Theta^{\text t}\end{array}\right)f
=\left(\begin{array}{cc}\U&\ov{\Theta}\end{array}\right)
\left(\begin{array}{c}\U^*f\\0\end{array}\right)
=\U\U^*f=\U g.
$$

Let us show that $g\in H^2(\C^r)$. Let $\g$ be a vector in $\C^r$. Since 
$\U^{\text t}$ is outer, it follows that there exists a sequence $\{\f_j\}$
of functions in $H^2(\C^n)$ such that $\{\U^{\text t}\f_j\}$ converges in $H^2(\C^r)$
to the function identically equal to $\g$, and so the sequence
$\{\f_j^{\text t}\U\}$ converges to the function identically equal to $\g^{\text t}$.
Hence, 
$$
\lim_{j\to\be}\{\f_j^{\text t}f\}=\lim_{j\to\be}\{\f_j^{\text t}\U g\}
=\g^{\text t}g
$$
in $L^1$, and so $\g^{\text t}g\in H^2$ for any constant vector $\g$. Consequently,
$g\in H^2(\C^r)$.

We have
$$
\Phi f=
\W^*\left(\begin{array}{c}U\U^*\U g\\0\end{array}\right)=
\left(\begin{array}{cc}\ov{\O}&\Xi\end{array}\right)
\left(\begin{array}{c}Ug\\0\end{array}\right)
=\left(\begin{array}{c}\ov{\O}Ug\\0\end{array}\right).
$$
Clearly, $f$ is a maximizing vector of $H_\Phi$ if and only if 
$\ov{\O}Ug\in H^2_-(\C^r)$ which is equivalent to the condition 
$\bar z\O\ov{Ug}\in H^2(\C^r)$. Since $\O^{\text t}$ is outer, we can apply
the same reasoning as above to show that $\bar z\ov{Ug}\in H^2(\C^r)$ which
is equivalent to the fact that $Ug\in H^2_-(\C^r)$. But the latter just
means that $g\in \Ker T_U$, and so $g$ is a maximizing vector of $H_U$. $\bl$

Suppose now that the matrix function $\Psi$ in the factorization \rf{4.6} also
satisfies the condition $\|H_\Psi\|_{\text e}<\|H_\Psi\|$. Then we can continue
this process, find a best analytic approximation $G$ of $\Psi$ and factorize
$\Psi-G$ as in \rf{4.6}. If we are able to continue this diagonalization process till
the very end, we construct the unique superoptimal approximation $Q$ of $\Phi$
and obtain a canonical factorization of $\Phi-Q$. 

Therefore we need an estimate of
$\|H_\Psi\|_{\text e}$. In [PT] in the case $r=1$ it was shown that 
$\|H_\Psi\|_{\text e}\le\|H_\Phi\|_{\text e}$. We want to obtain the same inequality
for an arbitrary $r$. We could try to generalize the proof
given in [PT] to the case of an arbitrary $r$. However, we are going to choose 
another way. We would like to deduce this result for an arbitrary $r$ from the 
corresponding result in the case $r=1$.

\

\setcounter{equation}{0}
\setcounter{section}{4}
\section{\bf Relations with thematic factorizations}

\

In this section we compare partial canonical factorizations obtained in \S 4 with
partial thematic factorizations and find useful relations between the complementing
matrix functions $\Theta$ and $\Xi$ in \rf{4.5} and the corresponding complementing 
matrix functions in partial thematic factorizations.

Let $\Phi$ be a matrix function in $L^\be(\mm_{m,n})$ such that 
$\|H_\Phi\|_{\text e}<\|H_\Phi\|$. Let $r$ be the number of superoptimal
singular values of $\Phi$ equal to $\|H_\Phi\|$. Suppose that $r<\min\{m,n\}$.
It follows from the results of \S 4 that if $F$ is a best approximation of $\Phi$
by bounded analytic matrix functions, then $\Phi-F$ admits a partial canonical 
factorization
\beq
\label{5.1}
\Phi-F=\W^*\left(\begin{array}{cc}t_0U&0\\0&\Psi\end{array}\right)\V^*,
\end{equation}
where $\V$ and $\W^{\text t}$ are $r$-balanced matrix functions of the form \rf{4.5}
and \linebreak$\|\Psi\|_{L^\be}\le t_0=\|H_\Phi\|$, and $\|H_\Psi\|<t_0$. 

On the other hand, it follows from the results of [PT] that $\Phi-F$ admits a
partial thematic factorization of the form
\beq
\label{5.2}
\Phi-F=W_0^*\cdots W_{r-1}^*\left(\begin{array}{ccccc}
t_0u_0&0&\cdots&0&0\\0&t_0u_1&\cdots&0&0\\\vdots&\vdots&\ddots&\vdots&\vdots\\
0&0&\cdots&t_0u_{r-1}&0\\0&0&\cdots&0&\D\end{array}\right)
V^*_{r-1}\cdots V_0^*,
\end{equation}
where $\|\D\|_{L^\be}\le t_0$ and $\|H_\D\|<t_0$,
$$
V_j=\left(\begin{array}{cc}I_j&0\\0&\breve{V}_j\end{array}\right),\quad
W_j=\left(\begin{array}{cc}I_j&0\\0&\breve{W}_j\end{array}\right),\quad
0\le j\le r-1,
$$
with thematic matrix functions $\breve{V}_j,\breve{W}_j^{\text t}$,
and the $u_j$ are unimodular functions such that the Toeplitz operators
$T_{u_j}$ are Fredholm with $\ind T_{u_j}>0$. For $j=0$ we assume that
$\breve{V}_0=V_0$ and $\breve{W}_0=W_0$.

Suppose that
\beq
\label{5.a}
\breve{V}_j=\left(\begin{array}{cc}\bs{v}_j&\ov{\Theta}_j\end{array}\right),\quad
\breve{W}^{\text t}_j=\left(\begin{array}{cc}\bs{w}_j&\ov{\Xi}_j\end{array}\right),
\quad0\le j\le r-1,
\end{equation}
and 
$$
\V=\left(\begin{array}{cc}\U&\ov{\Theta}\end{array}\right),\quad
\W^{\text t}=\left(\begin{array}{cc}\O&\ov{\Xi}\end{array}\right),
$$
where the matrix functions $\bs{v}_j,\Theta_j,\bs{w}_j,\Xi_j, \U,\Theta,\O,\Xi$
are inner and co-outer.

\begin{thm}
\label{t5.1}
Under the above hypotheses there exist constant unitary matrices 
$\frak U_1\in\mm_{n-r,n-r}$ and $\frak U_2\in\mm_{m-r,m-r}$ such that
\beq
\label{5.3}
\Theta=\Theta_0\Theta_1\cdots\Theta_{r-1}\frak U_1
\end{equation}
and
\beq
\label{5.4}
\Xi=\Xi_0\Xi_1\cdots\Xi_{r-1}\frak U_2.
\end{equation}
\end{thm}

\Pf It follows from Theorem \ref{t3.4} that if we replace $F$ with another best 
approximation $G$, the matrix function $\Phi-G$ will still admit factorizations of 
the forms \rf{5.1} and \rf{5.2} with the same matrix functions $\V,\W,V_j,W_j$.
Hence, we may assume that $F\in\bs{\O}_r$ (see \S 1). Then the matrix functions
$\Psi$ in \rf{5.1} and $\D$ in \rf{5.2} satisfy the inequalities
$$
\|\Psi\|_{L^\be}<t_0\quad\mbox{and}\quad\|\D\|_{L^\be}<t_0.
$$

Define the function $\r:\R\to\R$ by
$$
\r(t)=\left\{\begin{array}{ll}t,&t\ge t_0^2\\0,&t<t_0^2\end{array}\right..
$$
Consider the operator $M:H^2(\C^n)\to L^2(\C^n)$ of multiplication by the
matrix function $\r\big((\Phi-F)^{\text t}(\ov{\Phi-F})\big)$. It was proved in 
Lemmas 3.4 and 3.5 of [AP1] that
\beq
\label{5.5}
\Ker M=\Theta_0\Theta_1\cdots\Theta_{r-1}H^2(\C^{n-r}).
\end{equation}
On the other hand, it follows from \rf{5.1} that
$$
(\Phi-F)^{\text t}(\ov{\Phi-F})=\ov{\V}
\left(\begin{array}{cc}t^2_0I_r&0\\0&\Psi^{\text t}\ov{\Psi}\end{array}\right)
\V^{\text t}
$$
and since $\|\Psi^{\text t}\ov\Psi\|_{L^\be}<t_0^2$, we have
$$
\r\big((\Phi-F)^{\text t}(\ov{\Phi-F})\big)=\ov{\V}
\left(\begin{array}{cc}t^2_0I_r&0\\0&0\end{array}\right)
\V^{\text t}=\ov{\V}
\left(\begin{array}{cc}t^2_0I_r&0\\0&0\end{array}\right)
\left(\begin{array}{c}\U^{\text t}\\\Theta^*\end{array}\right).
$$
By Theorem \ref{t3.1}, $\Ker T_{\U^{\text t}}=\Theta H^2(\C^{n-r})$. It is easy to
see now that \linebreak$\Ker M=\Theta H^2(\C^{n-r})$. Together with \rf{5.5} this 
yields
$$
\Theta_0\Theta_1\cdots\Theta_{r-1}H^2(\C^{n-r})=\Theta H^2(\C^{n-r})
$$
which means that both inner functions $\Theta$ and 
$\Theta_0\Theta_1\cdots\Theta_{r-1}$ determine the same invariant 
subspace of multiplication by $z$ on $H^2(\C^n)$. Therefore there exists a constant 
unitary function $\frak U_1$ such that \rf{5.3} holds (see \S 2). To prove \rf{5.4}
we can apply \rf{5.3} to $(\Phi-F)^{\text t}$. $\bl$

\begin{cor}
\label{t5.2}
Let $\Psi$ and $\D$ be the matrix functions in the factorizations {\em\rf{5.1}}
and {\em\rf{5.2}}. Then
\beq
\label{5.6}
\D={\frak U}_2\Psi{\frak U}_1^{\text t},
\end{equation}
where $\frak U_1$ and $\frak U_2$ are unitary matrices from {\em\rf{5.3}} and 
{\em\rf{5.4}}.
\end{cor}

\Pf By Corollary 3.2 of [AP1],
$$
\D=\Xi_{r-1}^*\cdots\Xi_1^*\Xi_0^*(\Phi-F)\ov{\Theta_0\Theta_1\cdots\Theta_{r-1}}.
$$
By Theorem \ref{t5.1},
$$
\D={\frak U}_2\Xi^*(\Phi-F)\ov{\Theta}{\frak U}_1^{\text t}.
$$
On the other hand, it is easy to see from \rf{5.1} that
\beq
\label{5.b}
\Psi=\Xi^*(\Phi-F)\ov{\Theta},
\end{equation}
which implies \rf{5.6}. $\bl$

Now we are in a position to estimate $\|H_\Psi\|_{\text e}$ for the 
residual entry $\Psi$ in the factorization \rf{5.1}. This will be used in
\S 7 to obtain canonical factorizations of very badly approximable matrix functions.

\begin{thm}
\label{t5.3}
Let $\Phi$ be a function in $L^\be(\mm_{m,n})$ such that 
{\em$\|H_\Phi\|_{\text e}<\|H_\Phi\|$}. Then for the residual entry $\Psi$ in the
partial canonical factorization {\em\rf{5.1}} the following inequality holds {\em
$$
\|H_\Psi\|_{\text e}\le\|H_\Phi\|_{\text e}.
$$}
\end{thm}

\Pf Iterating Theorem 6.3 of [PT], we find that 
$\|H_\D\|_{\text e}\le\|H_\Phi\|_{\text e}$.
The result follows now from \rf{5.6}. $\bl$

Consider now the unitary-valued matrix function $U$ in the partial canonical
factorization \rf{5.1}. As we have observed in the remark following the proof
of Theorem \ref{t4.3}, the Toeplitz operator $T_U$ is Fredholm. We are going to 
evaluate now the index of $T_U$ in terms of the indices of the partial thematic 
factorization \rf{5.2}. Recall that the indices $k_j$ of the partial thematic
factorization \rf{5.2} are defined by
$$
k_j\df\ind T_{u_j},\quad 0\le j\le r-1.
$$

\begin{thm}
\label{t5.4}
Let $\Phi\in L^\be(\mm_{m,n})$ and {\em$\|H_\Phi\|_{\text e}<\|H_\Phi\|$}.
Then the entry $U$ of of the diagonal block matrix function in
the partial canonical factorization {\em\rf{5.1}} satisfies
$$
\ind T_U=\dim\Ker T_U=k_0+k_1+\cdots+k_{r-1}.
$$
\end{thm}

\Pf Without loss of generality we may assume that $\|H_\Phi\|=1$.
By Theorem \ref{t4.3}, $U$ is very badly approximable, $\|H_U\|_{\text e}<\|H_U\|=1$.
Therefore the Toeplitz operator $T_{zU}$ has dense range in $H^2(\C^r)$ (see
Theorem 3.1 of [PY1] where this fact was proved in the case 
$U\in(H^\be+C)(\mm_{r,r})$, however the proof given in [PY1] also works
in the more general case $\|H_U\|_{\text e}<\|H_U\|=1$). Therefore
$\Ker T^*_U=\{0\}$, and so $\ind T_U=\dim\Ker T_U$.

By Theorem 9.3 of [PT],
\beq
\label{5.8}
\dim\{f\in H^2(\C^n):~\|H_\Phi f\|_2=\|f\|_2\}=k_0+k_1+\cdots+k_{r-1}
\end{equation}
(earlier this result was proved in [PY2], Theorem 2.2 in the case 
\linebreak$\Phi\in(H^\be+C)(\mm_{m,n})$). Let us show that the left-hand side
of \rf{5.8} is equal to $\dim\Ker T_U$. Indeed, if $g\in H^2(\C^r)$, then
$g\in\Ker T_U$ if and only if $g$ is a maximizing vector of $H_U$.
By Lemma \ref{t4.8}, 
$$
\dim\{g\in H^2(\C^r):~\|H_Ug\|_2=\|g\|_2\}=
\dim\{f\in H^2(\C^n):~\|H_\Phi f\|_2=\|f\|_2\}
$$
which proves the result. $\bl$

To conclude this section, we obtain an inequality between the singular values
of $H_\Phi$ and the singular values of $H_\Psi$ in the factorization \rf{5.1}.

\begin{thm}
\label{t5.5}
Let $\Phi\in L^\be(\mm_{m,n})$ and {\em$\|H_\Phi\|_{\text e}<\|H_\Phi\|$}. Suppose 
that \linebreak$F\in H^\be(\mm_{m,n})$ is a best approximation of $\Phi$ by bounded 
analytic functions and $\Phi-F$ is represented by the partial canonical factorization
{\em\rf{5.1}}. Then
\beq
\label{5.9}
s_j(H_\Psi)\le s_{\iota+j}(H_\Phi),\quad j\ge0,
\end{equation}
where $\iota\df\ind T_U$.
\end{thm}

Note that \rf{5.9} was proved in [PT], Theorem 10.1 in the case $r=1$.
We reduce the general case to the case $r=1$.

\Pf Consider the partial thematic factorization \rf{5.2}. Let $k_j\df\ind T_{u_j}$,
$0\le j\le r-1$, be the indices of this factorization. We have
\beq
\label{5.10}
\Phi-F=W_0^*\left(\begin{array}{cc}t_0u_0&0\\0&\Phi^{[1]}\end{array}\right)V_0^*,
\end{equation}
where $\Phi^{[1]}$ is given by the partial thematic factorization
$$
\Phi^{[1]}=\breve W_1^*\cdots 
\left(\begin{array}{cc}I_{r-2}&0\\0&\breve W_{r-1}^*\end{array}\right)
\left(\begin{array}{cccc}
t_0u_1&\cdots&0&0\\\vdots&\ddots&\vdots&\vdots\\
0&\cdots&t_0u_{r-1}&0\\0&\cdots&0&\D\end{array}\right)
\left(\begin{array}{cc}I_{r-2}&0\\0&\breve V_{r-1}^*\end{array}\right)
\cdots\breve V_1^*.
$$
Now we can apply Theorem 10.1 of [PT] to the factorization \rf{5.10} and find that
$$
s_j(H_{\Phi^{[1]}})\le s_{k_0+j}(H_\Phi),\quad j\ge0.
$$
Then we can apply Theorem 10.1 of [PT] to the above partial thematic factorization 
of $\Phi^{[1]}$, {\it etc}. After applying Theorem 10.1 of [PT]
$r$ times we obtain the inequality
$$
s_j(H_\D)\le s_{k_0+\cdots+k_{r-1}+j}(H_\Phi),\quad j\ge0.
$$
The result follows now from Corollary \ref{t5.2} and Theorem \ref{t5.4}. $\bl$

\

\setcounter{equation}{0}
\setcounter{section}{5}
\section{\bf Invertibility of the Toeplitz operators $\bs{T_\V}$ and $\bs{T_\W}$}

\

For a matrix function $\Phi\in L^\be(\mm_{m,n})$ satisfying
$\|H_\Phi\|_{\text e}<\|H_\Phi\|$ we have constructed in \S 4 balanced matrix 
functions $\V$ and $\W^{\text t}$. We have mentioned in \S 3 that by Lemma 6.2 of
[Pe4], $T_\V$ and $T_\W$ have trivial kernels and dense ranges. In this section
we prove that under the condition $\|H_\Phi\|_{\text e}<\|H_\Phi\|$ the Toeplitz
operators $T_\V$ and $T_\W$ are invertible. In the case $r=1$ it was proved in
[PT] (see the proof of Theorem 5.1 of [PT]) that the Toeplitz operators $T_\V$ and 
$T_{\W^{\text t}}$ are invertible. Instead of trying to generalize the proof given
in [PT] to the case of an arbitrary $r$ we are going to reduce the general case to 
the case $r=1$.

Then we are going to prove that the matrix functions $\U,\Theta,\O,\Xi$ given by
\rf{4.5} are left invertible in $H^\be$. Note that this result in the case $r=1$
plays an important role (see [PY2] and [PT]).

We are going to use the following fact from [Pe2] (see also [PK] where the
scalar case was considered):

{\it If $V$ is an $n\times n$ unitary-valued function such that the Toeplitz 
operator \linebreak$T_V:H^2(\C^n)\to H^2(\C^n)$
has trivial kernel and dense range, then the operators $H_{V^*}^*H_{V^*}$
and $H^*_VH_V$ are unitarily equivalent. In particular, $\|H_V\|=\|H_{V^*}\|$.}

We need another well known fact (see [D] or [Ni] for the scalar case, in the
matrix case the proof is the same):

{\it If $V$ is an $n\times n$ unitary-valued function, then $T_V$ is invertible
on $H^2(\C^n)$ if and only if $\|H_V\|<1$ and $\|H_{V^*}\|<1$.}

Consider now a balanced matrix function 
$\V=\left(\begin{array}{cc}\U&\ov{\Theta}\end{array}\right)$. Clearly,
$\|H_\V\|=\|H_{\ov\Theta}\|$ and $\|H_{\V^*}\|=\|H_{\ov\U}\|$. Since $T_\V$ has
trivial kernel and dense range (Lemma 6.2 of [Pe4]), it follows that 
\beq
\label{6.1}
\|H_{\ov\Theta}\|=\|H_{\ov\U}\|.
\end{equation}
$T_\V$ is invertible if and only if the norms in \rf{6.1} are less than 1. It is also 
clear that the last condition is equivalent to the invertibility of 
$T_{\V^{\text t}}$.

\begin{thm}
\label{t6.1}
Let $\Phi\in L^\be(\mm_{m,n})$ and suppose that 
{\em$\|H_\Phi\|_{\text e}<\|H_\Phi\|$}. Let $\V$ and $\W^{\text t}$ be the 
$r$-balanced matrix functions in the partial canonical factorization {\em\rf{5.1}}.
Then the Toeplitz operators $T_\V$, {\em$T_{\V^{\text t}}$}, $T_{\W}$, 
and {\em$T_{\W^{\text t}}$} are invertible.
\end{thm}

Recall that for $r=1$ this was proved in [PT] (see the proof of Theorem 5.1 of [PT]).
We are going to reduce the general case to the case $r=1$.

\Pf Consider the partial thematic factorization \rf{5.2} and the corresponding
thematic matrix functions $\breve{V}_j$ and the inner matrix functions $\Theta_j$ 
given by \rf{5.a}. As we have already mentioned, it was shown in the proof of 
Theorem 5.1 of [PT] that the Toeplitz operators $T_{\breve{V}_j}$ are invertible for 
$0\le j\le r-1$. Since the $\breve V_j$ are 1-balanced, this is equivalent to the 
fact that $\|H_{\ov\Theta_j}\|<1$, $0\le j\le r-1$ (see the discussion preceding the
statement of Theorem \ref{t6.1}). To prove that $T_\V$ is invertible, it is 
sufficient to show that $\|H_{\ov\Theta}\|<1$, where as usual, 
$\V=\left(\begin{array}{cc}\U&\ov{\Theta}\end{array}\right)$. 

By Theorem \ref{t5.1}, $\Theta=\Theta_0\cdots\Theta_{r-1}\frak U$, where
$\frak U$ is a constant unitary matrix. Clearly, to show that $\|H_{\ov\Theta}\|<1$,
it is sufficient to prove the following lemma.

\begin{lem}
\label{t6.2}
Let $\varTheta_1$ be an $n\times k$ inner matrix function and let $\varTheta_2$ be 
a $k\times l$ inner matrix function such that $\|H_{{\ov\varTheta}_1}\|<1$
and $\|H_{{\ov\varTheta}_2}\|<1$. Let $\varTheta=\varTheta_1\varTheta_2$. Then
$\|H_{\ov\varTheta}\|<1$.
\end{lem}

Let us first complete the proof of Theorem \ref{t6.1}. Applying Lemma
\ref{t6.2} inductively, we  find that $\|H_{\ov\Theta}\|<1$. As we have explained
this is equivalent to the invertibility of $T_\V$ and $T_{\V^{\text t}}$. Similarly,
one can prove that the Toeplitz operators $T_{\W^{\text t}}$ and $T_\W$ are 
invertible. $\bl$

{\bf Proof of Lemma \ref{t6.2}.} Let $\s<1$ be a positive number such that
$\|H_{{\ov\varTheta}_1}\|<\s$ and $\|H_{{\ov\varTheta}_2}\|<\s$. Let 
$f\in H^2(\C^l)$. We have
\begin{eqnarray*}
\|H_{\ov\varTheta}f\|^2_2&=&\|\pp_-\ov{\varTheta_1}\ov{\varTheta_2}f\|_2^2\\
&=&\|\pp_-\ov{\varTheta_1}\pp_-\ov{\varTheta_2}f+
\pp_-\ov{\varTheta_1}\pp_+\ov{\varTheta_2}f\|_2^2\\
&=&\|\ov{\varTheta_1}\pp_-\ov{\varTheta_2}f+
\pp_-\ov{\varTheta_1}\pp_+\ov{\varTheta_2}f\|_2^2\\
&=&\|\pp_-\ov{\varTheta_2}f+
\varTheta_1^{\text t}\pp_-\ov{\varTheta_1}\pp_+\ov{\varTheta_2}f\|_2^2.
\end{eqnarray*}
We claim that the functions $\pp_-\ov{\varTheta_2}f$ and 
$\varTheta_1^{\text t}\pp_-\ov{\varTheta_1}\pp_+\ov{\varTheta_2}f$ are orthogonal.
Indeed, let $\f\in H^2_-(\C^k)$ and $\psi\in H^2(\C^k)$. We have
$$
(\f,\varTheta_1^{\text t}\pp_-\ov{\varTheta_1}\psi)=
(\ov{\varTheta_1}\f,\pp_-\ov{\varTheta_1}\psi)=
(\ov{\varTheta_1}\f,\ov{\varTheta_1}\psi)=(\f,\psi)=0, 
$$ 
since $\ov{\varTheta_1}$ takes isometric values almost everywhere on $\T$. It follows
that
\begin{eqnarray*}
\|H_{\ov\varTheta}f\|^2_2&=&\|\pp_-\ov{\varTheta_2}f\|^2_2+
\|\varTheta_1^{\text t}\pp_-\ov{\varTheta_1}\pp_+\ov{\varTheta_2}f\|_2^2\\
&\le&\|\pp_-\ov{\varTheta_2}f\|^2_2+
\|\pp_-\ov{\varTheta_1}\pp_+\ov{\varTheta_2}f\|_2^2\\
&=&\|\pp_-{\ov\varTheta_2}f\|^2_2+
\|H_{\ov{\varTheta_1}}\pp_+\ov{\varTheta_2}f\|_2^2\\
&\le&\|\pp_-{\ov\varTheta_2}f\|^2_2+\s^2\|\pp_+{\ov\varTheta_2}f\|^2_2\\
&=&\s^2(\|\pp_-{\ov\varTheta_2}f\|^2_2+\|\pp_+{\ov\varTheta_2}f\|^2_2)
+(1-\s^2)\|\pp_-{\ov\varTheta_2}f\|^2_2\\
&=&\s^2\|\ov{\varTheta_2}f\|^2_2+(1-\s^2)\|H_{\ov\varTheta_2}f\|^2_2\\
&\le&\s^2\|f\|^2_2+\s^2(1-\s^2)\|f\|^2_2=(2\s^2-\s^4)\|f\|^2_2.
\end{eqnarray*}
The result follows now from the trivial inequality $2\s^2-\s^4<1$. $\bl$

Now we are going to prove that the matrix functions $\U,\Theta,\O,\Xi$
are left invertible in $H^\be$.
Recall that a matrix function $\D\in H^\be(\mm_{\iota,\kappa})$ is called
{\it left invertible in} $H^\be$ if there exists a matrix function 
$\L\in H^\be(\mm_{\kappa,\iota})$ such that $\L\D$ is identically equal to 
$I_\kappa$. Recall that by a theorem of Sz.-Nagy and Foias [SNF2],
$\D$ is left invertible in $H^\be$ if and only if the Toeplitz operator
$T_{\ov\D}:H^2(\C^\kappa)\to H^2(\C^\iota)$ is bounded below, i.e.,
$\|T_{\ov\D}\f\|_2\ge\d\|\f\|_2$ for some $\d>0$.

\begin{thm}
\label{t6.3}
Let $\Phi\in L^\be(\mm_{m,n})$ and {\em$\|H_\Phi\|_{\text e}<\|H_\Phi\|$}.
Suppose that $\V$ and $\W$ are matrix functions of the form {\em\rf{4.5}}
in the partial canonical factorization {\em\rf{5.1}}. 
Then $\U,\Theta,\O,\Xi$ are left invertible in $H^\be$.
\end{thm}

\Pf Let us prove that $\Theta$ is left invertible. By Theorem \ref{t6.1}.
the Toeplitz operator $T_\V$ is invertible. Hence, $\|T_\V f\|_2\ge\d\|f\|_2$, 
$f\in H^2(\C^n)$, for some $\d>0$. If we apply this inequality to functions $f$ of 
the form $f=\left(\begin{array}{c}0\\\f\end{array}\right)$, we find that
$\|T_{\ov{\Theta}}\f\|_2\ge\d\|\f\|_2$, $\f\in H^2(\C^{n-r})$. It follows now from
the Sz.-Nagy--Foias theorem mentioned above that $\Theta$ is left invertible in
$H^\be$. 

To prove that $\U$ is left invertible in $H^\be$, we can apply the above reasoning to
the matrix function $\ov\V$. Finally, to obtain the same results for $\O$ and $\Xi$,
we can apply the above reasoning to the matrix functions $\W^*$ and $\W^{\text t}$.
$\bl$

\

\setcounter{equation}{0}
\setcounter{section}{6}
\section{\bf Canonical factorizations of very badly approximable functions}

\

Given $\Phi\in L^\be(\mm_{m,n})$, consider the sequence $\{t_j\}$ of its
superoptimal singular values. Suppose that
\beq
\label{7.1}
t_0=\cdots=t_{r_1-1}>t_{r_1}=\cdots=t_{r_2-1}>\cdots>
t_{r_{\iota-1}}=\cdots=t_{r_\iota-1}
\end{equation}
are all nonzero superoptimal singular values of $\Phi$, i.e., $t_0$ has
multiplicity $r_1$ and $t_{r_j}$ has multiplicity $r_{j+1}-r_j$, $1\le j\le\iota-1$.

By Theorem \ref{t4.3}, if $\|H_\Phi\|_{\text e}<\|H_\Phi\|$ and $F_0$ is a best 
approximation of $\Phi$ by bounded analytic matrix function, then $\Phi-F_0$
admits a factorization
$$
\Phi-F_0=\W_0^*\left(\begin{array}{cc}t_0U_0&0\\0&\Phi^{[1]}\end{array}\right)\V_0^*,
$$
where $\V_0$ and $\W^{\text t}_0$ are $r_1$-balanced matrix functions,
$U_0$ is very badly approximable unitary-valued matrix function such that
$\|H_{U_0}\|_{\text e}<1$, and $\Phi^{[1]}$ is a matrix function in 
$L^\be(\mm_{m-r_1,n-r_1})$ such that $\|H_{\Phi^{[1]}}\|=t_{r_1}$.

By Theorem \ref{t5.3}, $\|H_{\Phi^{[1]}}\|_{\text e}\le\|H_\Phi\|_{\text e}$. If
$t_{r_1}$ is still greater than $\|H_\Phi\|_{\text e}$, we can apply Theorem
\ref{t4.3} to $\Phi^{[1]}$ and find that for a best approximation $G_1$ of
$\Phi^{[1]}$ the matrix function $\Phi^{[1]}-G_1$ admits a factorization
$$
\Phi^{[1]}-G_1=\breve{\W}_1^*
\left(\begin{array}{cc}t_{r_1}U_1&0\\0&\Phi^{[2]}\end{array}\right)
\breve{\V}_1^*,
$$
where $\breve{\V}_1$ and $\breve{\W}^{\text t}_1$ are $(r_2-r_1)$-balanced matrix 
functions, $U_1$ is a very badly approximable unitary-valued matrix function 
of size $(r_2-r_1)\times(r_2-r_1)$ such that
$\|H_{U_1}\|_{\text e}<1$, and $\Phi^{[2]}$ is a matrix function in 
$L^\be(\mm_{m-r_2,n-r_2})$ such that \linebreak$\|H_{\Phi^{[2]}}\|=t_{r_2}$.

We can apply now Theorem \ref{t3.4} and find a matrix function 
$F_1\in H^\be(\mm_{m,n})$ such that
$$
\Phi-F_1=\W_0^*\W_1^*
\left(\begin{array}{ccc}t_0U_0&0&0\\0&t_{r_1}U_1&0\\
0&0&\Phi^{[2]}\end{array}\right)
\V_1^*\V_0^*,
$$
where
$$
\V_1=\left(\begin{array}{cc}I_{r_1}&0\\0&\breve{\V}_1\end{array}\right)
\quad\mbox{and}\quad
\W_1=\left(\begin{array}{cc}I_{r_1}&0\\0&\breve{\W}_1\end{array}\right).
$$

If $t_{r_2}$ is still greater than $\|H_\Phi\|_{\text e}$, we can continue this
process and apply Theorem \ref{t4.3} to $\Phi^{[2]}$. Suppose now that 
$t_{r_{d-1}}>\|H_\Phi\|_{\text e}$, $2\le d\le \iota$. Then continuing the above 
process and applying Theorem \ref{t3.4}, we can find a function 
$F\in H^\be(\mm_{m,n})$ such that $\Phi-F$ admits a factorization
\beq
\label{7.2}
\Phi-F=\W_0^*\cdots\W^*_{d-1}
\left(\begin{array}{ccccc}t_0U_0&0&\cdots&0&0\\0&t_{r_1}U_1&\cdots&0&0\\
\vdots&\vdots&\ddots&\vdots&\vdots\\
0&0&\cdots&t_{r_{d-1}}U_{d-1}&0\\
0&0&\cdots&0&\Phi^{[d]}\end{array}\right)
\V^*_{d-1}\cdots\V_0^*,
\end{equation}
where the $U_j$ are $(r_{j+1}-r_j)\times(r_{j+1}-r_j)$ very badly approximable 
unitary valued functions such that $\|H_{U_j}\|_{\text e}<1$,
\beq
\label{7.3}
\V_j=\left(\begin{array}{cc}I_{r_j}&0\\0&\breve{\V}_j\end{array}\right)
\quad\mbox{and}\quad
\W_j=\left(\begin{array}{cc}I_{r_j}&0\\0&\breve{\W}_j\end{array}\right),
\quad1\le j\le d-1,
\end{equation}
$\breve{\V}_j$ and $\breve{\W}^{\text t}_j$ are $(r_{j+1}-r_j)$-balanced matrix 
functions, and $\Phi$ is a matrix function satisfying
\beq
\label{7.4}
\|\Phi^{[d]}\|_{L^\be}\le t_{r_{d-1}},\quad\mbox{and}\quad 
\|H_{\Phi^{[d]}}\|<t_{r_{d-1}}.
\end{equation}

Factorizations of the form \rf{7.2} with $\V_j$ and $\W_j$ of the form \rf{7.3} and
$\Phi^{[d]}$ satisfying \rf{7.4} are called {\it partial canonical factorizations}.
The matrix function $\Phi^{[d]}$ is called
the {\it residual entry} of the partial canonical factorization \rf{7.2}.

Finally, if $\|H_\Phi\|_{\text e}$ is less than the smallest nonzero 
superoptimal singular value $t_{r_\iota-1}$, then we can complete this process
and construct the unique superoptimal approximation of $\Phi$ by
bounded analytic matrix functions. This proves the following theorem.

\begin{thm}
\label{t7.1}
Let $\Phi\in L^\be(\mm_{m,n})$ and suppose that the nonzero superoptimal singular
values of $\Phi$ satisfy {\em\rf{7.1}}. If {\em$\|H_\Phi\|_{\text e}<t_{r_\iota-1}$}
and $F$ is the unique superoptimal 
approximation of $\Phi$ by bounded analytic functions, then $\Phi-F$ admits
a factorization
\beq
\label{7.5}
\Phi-F=\W_0^*\cdots\W^*_{\iota-1}
\left(\begin{array}{ccccc}t_0U_0&0&\cdots&0&0\\0&t_{r_1}U_1&\cdots&0&0\\
\vdots&\vdots&\ddots&\vdots&\vdots\\
0&0&\cdots&t_{r_{\iota-1}}U_{\iota-1}&0\\
0&0&\cdots&0&0\end{array}\right)
\V^*_{\iota-1}\cdots\V_0^*,
\end{equation}
where the $\V_j$ and $\W_j$, and $U_j$ are as above.
\end{thm}

Note that the lower right entry of the diagonal matrix function on the right-hand
side of \rf{7.5} is the zero matrix function of size $(m-r_\iota)\times(n-r_\iota)$.
Here it may happen that $m-r_\iota$ or $n-r_\iota$ can be zero.

Clearly, the left-hand side of \rf{7.5} is a very badly approximable matrix
function. Factorizations of the form \rf{7.5} are called {\it canonical 
factorizations} of very badly approximable matrix functions.

The following result shows that the right-hand side of \rf{7.5} is always 
a very badly approximable matrix function.

\begin{thm}
\label{t7.2}
Let $\Phi$ be a matrix function in $L^\be(\mm_{m,n})$ such that
\linebreak{\em$\|H_\Phi\|_{\text e}<\|H_\Phi\|$}, 
$\iota$ is a positive integer, $r_1,\cdots,r_\iota$ are positive integers satisfying
$$
r_1<r_2<\cdots<r_\iota
$$ 
and
$$
\s_0>\s_1>\cdots>\s_{\iota-1}>0.
$$
Suppose that $\Phi$ admits a factorization
$$
\Phi=\W_0^*\cdots\W^*_{\iota-1}
\left(\begin{array}{ccccc}\s_0U_0&0&\cdots&0&0\\0&\s_1U_1&\cdots&0&0\\
\vdots&\vdots&\ddots&\vdots&\vdots\\
0&0&\cdots&\s_{\iota-1}U_{\iota-1}&0\\
0&0&\cdots&0&0\end{array}\right)
\V^*_{\iota-1}\cdots\V_0^*,
$$
in which the $U_j$, $\V_j$, and $\W_j$ are as above.
Then $\Phi$ is very badly approximable and the superoptimal singular values of $\Phi$
are given by
$$
t_\kappa(\Phi)=\left\{\begin{array}{ll}\s_0,&\kappa<r_1\\
\s_j,&r_j\le \kappa<r_{j+1}
\\0,&\kappa\ge r_\iota
\end{array}\right..
$$
\end{thm}

\Pf Let
$$
\V_0=\left(\begin{array}{cc}\U&\ov{\Theta}\end{array}\right)\quad\mbox{and}\quad
\W_0=\left(\begin{array}{cc}\O&\ov{\Xi}\end{array}\right),
$$
where $\U,\Theta,\O,\Xi$ are inner and co-outer matrix functions.

By Theorem \ref{t4.7}, the minimal invariant subspace of multiplication by $z$
on $H^2(\C^n)$ that contains all maximizing vectors of $H_\Phi$ is $\U H^2(\C^{r_1})$
and the minimal invariant subspace of multiplication by $z$ on $H^2(\C^m)$ that
contains all maximizing vectors of $H_{\Phi^{\text t}}$ is $\O H^2(\C^{r_1})$.
By Theorem \ref{t4.6}, $\Phi$ is badly approximable and 
$$
t_0(\Phi)=\cdots=t_{r_1-1}=\s_0.
$$
By Theorem \ref{t3.4}, we can reduce our problem to the function
$$
\breve\W_1^*\cdots
\left(\begin{array}{cc}I_{r_{\iota-1}-r_1}&0\\0&\breve\W^*_{\iota-1}
\end{array}\right)\!\!
\left(\begin{array}{cccc}\s_1U_1&\cdots&\0&\0\\
\vdots&\ddots&\vdots&\vdots\\
\0&\cdots&\s_{\iota-1}U_{\iota-1}&\0\\
\0&\cdots&\0&\0\end{array}\right)\!\!
\left(\begin{array}{cc}I_{r_{\iota-1}-r_1}&0\\0&\breve\V^*_{\iota-1}
\end{array}\right)
\cdots\breve\V_1^*.
$$
This function is also represented by a canonical factorization which makes it
possible to continue this process and prove that $\Phi$ is very badly approximable
and the superoptimal singular values of $\Phi$ satisfy the desired equality. $\bl$

\

\setcounter{equation}{0}
\setcounter{section}{7}
\section{\bf Invariance properties of canonical factorizations}

\

In this section we demonstrate an advantage of canonical factorizations
over thematic factorizations. Namely, we show that canonical factorizations
possess certain invariance properties, i.e., the matrix functions $U_j$ 
in the canonical factorization \rf{7.5} are uniquely
determined modulo unitary constant  factors. Moreover, if 
$$
\V_0=\left(\begin{array}{cc}\U_0&\ov{\Theta}_0\end{array}\right),\quad
\W^{\text t}_0=\left(\begin{array}{cc}\O_0&\ov{\Xi}_0\end{array}\right),
$$
$$
\breve{\V}_j=\left(\begin{array}{cc}\U_j&\ov{\Theta}_j\end{array}\right),\quad
\breve{\W}^{\text t}_j=\left(\begin{array}{cc}\O_j&\ov{\Xi}_j\end{array}\right),
\quad1\le j\le\iota-1,
$$
and the $\breve{\V}_j$ and $\breve{\W}_j$ are given by \rf{7.3}, then
the matrix functions $\U_j,\Theta_j,\O_j,\Xi_j$ are also uniquely determined 
modulo constant unitary factors.

We start with partial canonical factorizations of the form \rf{5.1}. Suppose
that a matrix function $\Phi$ in $L^\be(\mm_{m,n})$ satisfies 
$\|H_\Phi\|_{\text e}<\|H_\Phi\|$ and admits partial canonical factorizations
\beq
\label{8.1}
\Phi=\left(\begin{array}{cc}\ov{\O}&\Xi\end{array}\right)
\left(\begin{array}{cc}\s U&0\\0&\Psi\end{array}\right)
\left(\begin{array}{cc}\U&\ov{\Theta}\end{array}\right)^*
\end{equation}
and
\beq
\label{8.2}
\Phi=\left(\begin{array}{cc}\ov{\O^\circ}&\Xi^\circ\end{array}\right)
\left(\begin{array}{cc}\s^\circ U^\circ&0\\0&\Psi^\circ\end{array}\right)
\left(\begin{array}{cc}\U^\circ&\ov{\Theta^\circ}\end{array}\right)^*,
\end{equation}
where $\|\Psi\|_{L^\be}\le\s$, $\|H_\Psi\|<\s$,
$\|\Psi^\circ\|_{L^\be}\le\s^\circ$, $\|H_{\Psi^\circ}\|<\s^\circ$,
$U$ is an $r\times r$
very badly approximable unitary-valued function such that $\|H_U\|_{\text e}<1$,
$U^\circ$ is an $r^\circ\times r^\circ$ very badly approximable unitary-valued 
function such that $\|H_{U^\circ}\|_{\text e}<1$, 
$\left(\begin{array}{cc}\U&\ov{\Theta}\end{array}\right)$ and
$\left(\begin{array}{cc}\O&\ov{\Xi}\end{array}\right)^{\text t}$ are
$r$-balanced matrix functions, and
$\left(\begin{array}{cc}\U^\circ&\ov{\Theta^\circ}\end{array}\right)$ and
$\left(\begin{array}{cc}\O^\circ&\ov{\Xi^\circ}\end{array}\right)^{\text t}$ are
$r^\circ$-balanced matrix functions.

\begin{thm}
\label{t8.1}
Let $\Phi$ be a badly approximable function in $L^\be(\mm_{m,n})$
such that {\em$\|H_\Phi\|_{\text e}<\|H_\Phi\|$}. Suppose that $\Phi$ admits
factorizations {\em\rf{8.1}} and {\em\rf{8.2}}. Then 
$r=r^\circ$, $\s=\s^\circ$ and there exist
unitary matrices $\frak V^\#,\frak V^\flat\in\mm_{r,r}$, $\frak U^\#\in\mm_{n-r,n-r}$,
$\frak U^\flat\in\mm_{m-r,m-r}$ such that
\beq
\label{8.3}
\U^\circ=\U\frak V^\#,\quad\O^\circ=\O\frak V^\flat,
\end{equation}
\beq
\label{8.4}
\Theta^\circ=\Theta\frak U^\#,\quad \Xi^\circ=\Xi\frak U^\flat,
\end{equation}
and
{\em
\beq
\label{8.5}
U^\circ=(\frak V^\flat)^{\text t}U\frak V^\#,\quad
\Psi^\circ=(\frak U^\flat)^*\Psi\ov{{\frak U}^\#}.
\end{equation}}
\end{thm}

\Pf By Theorem \ref{t4.6}, $\s=\s^\circ=\|H_\Phi\|$. Next, by Theorem \ref{t4.7},
the minimal invariant subspace of multiplication by $z$ on $H^2(\C^n)$ 
that contains all maximizing vectors of $H_\Phi$ is
equal to $\U H^2(\C^r)$ and at the same time it is equal to 
$\U^\circ H^2(\C^{r^\circ})$. It follows that $r=r^\circ$ and there exists
a unitary matrix $\frak V^\#\in\mm_{r,r}$ such that $\U^\circ=\U\frak V^\#$.
Applying the same reasoning to  $\Phi^{\text t}$, we find a unitary matrix function
$\frak V^\flat\in\mm_{r,r}$ such that $\O^\circ=\O\frak V^\flat$ which proves
\rf{8.3}.

By Theorem \ref{t3.1},
$$
\Theta H^2(\C^{n-r})=\Ker T_{\U^{\text t}}\quad\mbox{and}\quad
\Theta^\circ H^2(\C^{n-r})=\Ker T_{(\U^\circ)^{\text t}}.
$$
By \rf{8.3}, $\Ker T_{\U^{\text t}}=\Ker T_{(\U^\circ)^{\text t}}$ which implies
that there exist a unitary matrix $\frak U^\#\in\mm_{n-r,n-r}$ such that 
$\Theta^\circ=\Theta\frak U^\#$. Applying the same reasoning to $\Phi^{\text t}$,
we find a unitary matrix $\frak U^\flat\in\mm_{m-r,m-r}$ such that
$\Xi^\circ=\Xi\frak U^\flat$ which proves \rf{8.4}.

By \rf{4.7},
$$
\s U=\O^{\text t}\Phi\U\quad\mbox{and}\quad
\s U^\circ=(\O^\circ)^{\text t}\Phi\U^\circ.
$$
This implies the first equality in \rf{8.5}.

Finally, by \rf{5.b},
$$
\Psi=\Xi^*\Phi\ov{\Theta}\quad\mbox{and}\quad
\Psi^\circ={\Xi^\circ}^*\Phi\ov{\Theta^\circ}
$$
which completes the proof of \rf{8.5}. $\bl$

We can obtain similar results for arbitrary partial canonical factorizations
and canonical factorizations. Let us consider in detail the following special case.
Suppose that the matrix functions $\Psi$ and $\Psi^\circ$ in \rf{8.1} and \rf{8.2} 
admit the following partial canonical factorizations
\beq
\label{8.6}
\Psi=\left(\begin{array}{cc}\ov{\O}_1&\Xi_1\end{array}\right)
\left(\begin{array}{cc}\s_1 U_1&0\\0&\D\end{array}\right)
\left(\begin{array}{cc}\U_1&\ov{\Theta}_1\end{array}\right)^*
\end{equation}
and
\beq
\label{8.7}
\Psi^\circ=\left(\begin{array}{cc}\ov{\O_1^\circ}&\Xi_1^\circ\end{array}\right)
\left(\begin{array}{cc}\s_1^\circ U_1^\circ&0\\0&\D^\circ\end{array}\right)
\left(\begin{array}{cc}\U_1^\circ&\ov{\Theta_1^\circ}\end{array}\right)^*,
\end{equation}
where $\|\D\|_{L^\be}\le\s_1$, $\|H_\D\|<\s_1$
and $\|\D^\circ\|_{L^\be}<\s^\circ_1$, $\|H_{\D^\circ}\|<\s^\circ_1$.

By Theorem \ref{t8.1}, $\Psi^\circ=(\frak U^\flat)^*\Psi\ov{{\frak U}^\#}$. 
It follows now from \rf{8.7} that
\beq
\label{8.8}
\Psi=\left(\begin{array}{cc}\ov{\ov{\frak U^\flat}\O_1^\circ}&
\frak U^\flat\Xi_1^\circ\end{array}\right)
\left(\begin{array}{cc}\s_1^\circ U_1^\circ&0\\0&\D^\circ\end{array}\right)
\left(\begin{array}{cc}\ov{\frak U^\#}\U_1^\circ&
\ov{\frak U^\#\Theta_1^\circ}\end{array}\right)^*
\end{equation}
which is another partial canonical factorization of $\Psi$.

We can compare now the factorizations \rf{8.6} and \rf{8.8}. By Theorem
\ref{t8.1}, $\s_1^\circ=\s_1$ and there exist unitary matrices 
$\frak V_1^\#,\frak V_1^\flat,\frak U_1^\#,\frak U_1^\flat$ such that
$$
\U_1^\circ=(\frak U^\#)^{\text t}\U_1\frak V_1^\#,\quad
\O_1^\circ=(\frak U^\flat)^{\text t}\O_1{\frak V}_1^\flat,
$$
$$
\Theta_1^\circ=(\frak U^\#)^*\Theta_1\frak U_1^\#,\quad
\Xi_1^\circ=(\frak U^\flat)^*\Xi_1\frak U_1^\flat,
$$
and
$$
U_1^\circ=(\frak V_1^\flat)^{\text t}U_1\frak V_1^\#,\quad
\D^\circ=(\frak U_1^\flat)^*\D\ov{{\frak U_1}^\#}.
$$

It is easy to see that the same results hold in the case 
of arbitrary partial canonical factorizations as well as 
arbitrary canonical factorizations.

\

\setcounter{equation}{0}
\setcounter{section}{8}
\section{\bf Hereditary properties}

\

In this section we consider the following heredity problem. Suppose that the 
initial matrix function $\Phi$ belongs to a certain function space $X$. We study the
problem of whether all matrix functions in a (partial) canonical factorization
of $\Phi$ belong to the same space $X$ (by this we certainly mean that
all their entries are in $X$). This is not true for an arbitrary function
space $X$. In particular, it can be shown that this is not true if $X$ is the space
$C(\T)$ of continuous function. Nevertheless,
we prove that this is true for two natural classes of function spaces. 
The first class consists of so-called $\cR$-spaces, i.e., spaces that can be 
described in terms of rational approximation in the norm of $BMO$ (see [PK]). The 
second class of spaces consists of Banach algebras satisfying Axioms (A1)--(A4) 
below. In both cases we apply so-called recovery theorems for unitary-valued 
functions obtained in [Pe2] and [Pe3].

We are not going to give a formal definition of $\cR$-spaces and refer the reader
to [PK] for details. Roughly speaking, a linear $\cR$-space is a linear space $X$ 
of functions on $\T$ such that $X\subset VMO$ and there exists a K\"{o}the sequence 
space $E$ 
such that $\f\in X$ if and only if $\f\in BMO$ and the singular values
of the Hankel operators $H_\f$ and $H_{\bar\f}$ satisfy
$$
\{s_n(H_\f)\}_{n\ge0}\in E\quad\mbox{and}\quad\{s_n(H_{\bar\f})\}_{n\ge0}\in E.
$$
Recall that $E$ is a {\it K\"{o}the sequence space} if 
$$
\{c_n\}_{n\ge0}\in E,\quad|d_n|\le|c_n|\quad\Rightarrow\quad\{d_n\}_{n\ge0}\in E.
$$

Important examples of $\cR$-spaces are the Besov spaces $B_p^{1/p}$, $0<p<\be$ 
(see [PK], [Pe1] (the corresponding K\"{o}the space $E$ is $\ell^p$)
and the space $VMO$ of functions of vanishing mean oscillation ($E$ is the space 
$c_0$ of sequences converging to 0). We refer the reader to [G]  for definitions
and properties of the spaces $BMO$ and $VMO$.

The second class consists of function spaces $X\subset C(\T)$ 
that contain the trigonometric polynomials and satisfy the following axioms:

(A1) {\it If $f\in X$, then $\bar f\in X$ and $\pp_+f\in X$;}

(A2) {\it$X$ is a Banach algebra with respect to pointwise multiplication;}

(A3) {\it for every $\f\in X$ the Hankel operator $H_\f$ is a compact operator
from $X_+$ to $X_-$;}

(A4) {\it if $f\in X$ and $f$ does not vanish on $\T$, then $1/f\in X$.}

Here we use the notation
$$
X_+=\{f\in X:~\hat f(j)=0,~j<0\},\quad X_-=\{f\in X:~\hat f(j)=0,~j\ge0\}.
$$

For simplicity we write $\Phi\in X$ if all entries of a matrix function $\Phi$
belong to $X$.

Note here that the Besov classes $B_{p,q}^s$,
$1\le p<\be$, $1\le q\le\be$, $s>1/p$, the space of functions with absolutely
convergent Fourier series, the spaces
$$
\{f:~f^{(n)}\in VMO\},\quad \{f:~\pp_+f^{(n)}\in C(\T),~\pp_-f^{(n)}\in C(\T)\},
\quad n\ge1,
$$
satisfy (A1)--(A4).

Among nonseparable Banach spaces $X$ satisfying (A1)--(A4) we mention the 
following ones: the H\"{o}lder--Zygmund spaces $\L_\a$, $\a>0$, the spaces
$$
\{f:~f^{(n)}\in BMO\},\quad \{f:~\pp_+f^{(n)}\in H^\be,~\pp_-f^{(n)}\in H^\be\},
\quad n\ge1,
$$
the space
$$
\{f:~|\hat f(j)|\le\const(1+|j|)^{-\a}\},\quad\a>1,
$$
(see [Pe3]).

We need the following so-called recovery theorem for unitary-valued functions. 
Let $X$ be either an $\cR$-space or a space of functions 
satisfying (A1)--(A4) and let $\cal U$ be a unitary-valued matrix function in 
$X(\mm_{n,n})$ such that the Toeplitz operator $T_{\cal U}:H^2(\C^n)\to H^2(\C^n)$ 
has dense range. Then
\beq
\label{9.1}
\pp_-{\cal U}\in X\quad\Longrightarrow\quad{\cal U}\in X.
\end{equation}
Moreover, if $X$ is a Banach $\cR$-space, then we can conclude that
$$
\|\pp_-{\cal U}\|_X\le\const\|\cal U\|_X.
$$

For $\cR$-spaces this was proved in [Pe2], for spaces satisfying (A1)--(A4) this was
proved in [Pe3], see also [PK] for the scalar case. In fact, in [PK] and [Pe3] 
the assumptions on $X$ were slightly different but
it was shown in [AP2] that the method used in [Pe3] can be adjusted to work
for all spaces satisfying (A1)--(A4).

Note that both $\cR$-spaces and  spaces satisfying (A1)--(A4) are contained
in $VMO$. Therefore if $\Phi\in L^\be(\mm_{m,n})$ and $\pp_-\Phi\in X(\mm_{m,n})$, 
then $H_\Phi$ is compact. 

To establish the main result of this section we need the fact that
for a linear \linebreak$\cR$-space $X$ the space $X\cap L^\be$ forms an algebra
with respect to pointwise multiplication. Indeed, suppose that
$\f,\psi\in X\cap L^\be$. Clearly, $\f\in X$ if and only if $\bar\f\in X$, and so
it is sufficient to prove that $\pp_-(\f\psi)\in X$. 
Let $f\in H^2$. We have
\beq
\label{9.2}
H_{\f\psi}f=\pp_-\f\psi f=\pp_-\f\pp_+\psi f+\pp_-\f\pp_-\psi f=
H_\f T_\psi f+\breve{T}_\f H_\psi f,
\end{equation}
where the operator $\breve{T}_\f:H^2_-\to H^2_-$ is defined by 
$\breve{T}_\f g=\pp_-\f g$, $g\in H^2_-$. It is easy to deduce from \rf{9.2} that
$\pp_-(\f\psi)\in X$.

\begin{thm}
\label{t9.1}
Suppose that $X$ is a linear $\cR$-space or $X$ is a function space satisfying
{\em(A1)--(A4)}. Let $\Phi$ be a bounded $m\times n$ matrix 
function such that $\pp_-\Phi$ is a nonzero matrix function in $X(\mm_{m,n})$. 
If $F\in H^\be(\mm_{m,n})$ is a best approximation of $\Phi$ and $\Phi-F$ 
admits a partial canonical factorization
$$
\Phi-F=\left(\begin{array}{cc}\ov{\O}&\Xi\end{array}\right)
\left(\begin{array}{cc}t_0U&0\\0&\Psi\end{array}\right)
\left(\begin{array}{cc}\U&\ov{\Theta}\end{array}\right)^*,
$$
then $\U,\Theta,\O,\Xi,U,\pp_-\Psi\in X$.
\end{thm}

\Pf Assume without loss of generality that $t_0=1$. By Theorem \ref{t3.4}, if we 
replace a best approximating function $F$ with any other best approximation, we
do not change $\pp_-\Psi$. Thus we may assume that $F$ is the unique superoptimal
approximation of $\Phi$ by bounded analytic matrix functions. Then $\Phi-F$ belongs 
to $X$. For linear $\cR$-spaces this was proved in [PY1], Theorem 5.1. For spaces
$X$ satisfying (A1)--(A4) this is Theorem 9.2 of [Pe3]. (Theorem 9.2
is stated in [Pe3] under slightly different assumptions but using the results of \S 4
of [AP2], one can easily see that the proof given in [Pe3] works for all spaces 
satisfying (A1)--(A4).)

Let us first prove that $\Theta\in X$. Consider a partial thematic factorization
of $\Phi-F$ of the form \rf{5.2}. Then the matrix functions 
$V_j$ given by \rf{5.a} belong to $X$. If $X$ is a linear $\cR$-space, this
was proved in \S 5 of [PY1], for spaces satisfying (A1)--(A4) this was proved
in \S 9 of [Pe3]. (Again this was proved in [Pe3]
under slightly different assumptions but the results of \S 4 of [AP2] show that the 
proof given in [Pe3] works for all spaces satisfying (A1)--(A4).) In particular
the inner matrix functions $\Theta_j$ in \rf{5.a} belong to $X$. 
Since $X\cap L^\be$ is an algebra, it follows now
from \rf{5.3} that $\Theta\in X$.

Consider now the unitary-valued matrix function 
$\V=\left(\begin{array}{cc}\U&\ov{\Theta}\end{array}\right)$. We have just proved 
that $\pp_-\V\in X$. As we have mentioned in \S 3, the Toeplitz operator $T_\V$ has 
dense range in $H^2(\C^n)$. Therefore by \rf{9.1}, $\V\in X$, and so $\U\in X$.

If we apply the above reasoning to $\Phi^{\text t}$, we prove that $\Xi\in X$
and $\O\in X$. It follows now from \rf{4.7} that $U\in X$. Finally,
it follows from \rf{5.b} that $\pp_-\Psi\in X$. $\bl$

{\bf Remark.} It can be shown that if $X$ is a linear $\cR$-space, then the $X$-norms
of $\U,\Theta,\O,\Xi,U,\pp_-\Psi$ can be estimated in terms of the $X$-norm
of $\pp_-\Phi$.

Clearly, it follows from Theorem \ref{t9.1} that the same result holds for arbitrary
partial canonical factorizations. In particular, the following theorem holds.

\begin{thm}
\label{t9.2}
Let $\Phi\in L^\be(\mm_{m,n})$ and $F$ is the superoptimal approximation
of $\Phi$ by bounded analytic matrix functions. If {\em\rf{7.5}} 
is a canonical factorization of $\Phi-F$, then all factors on the right-hand side
of {\em\rf{7.5}} belong to $X$.
\end{thm}

Consider now separately the important case $X=VMO$. As we have already noted, $VMO$ 
is an $\cR$-space. It is well known that if $\Phi\in L^\be$, then $\pp_-\Phi\in VMO$
if and only if $\Phi\in H^\be+C$. It is also well known that $QC=VMO\cap L^\be$.

\begin{thm}
\label{t9.3}
Let $\Phi\in(H^\be+C)(\mm_{m,n})$ and $\pp_-\Phi\neq 0$. If $F$ is a best 
approximation of $\Phi$ by bounded analytic functions and $\Phi-F$ admits 
a partial canonical factorization {\em\rf{4.6}}, then $\V,\W,U\in QC$
and $\Psi\in H^\be+C$.
\end{thm}

Theorem \ref{t9.3} follows immediately from Theorem \ref{t9.1} if we put $X=VMO$.

\

\setcounter{equation}{0}
\setcounter{section}{9}
\section{\bf Very badly approximable unitary-valued functions}

\

In this section we study unitary-valued very badly approximable matrix functions
which play a crucial role in canonical factorizations. The following theorem
describes such functions $U$ under the condition $\|H_U\|_{\text e}<1$.

\begin{thm}
\label{10.1}
Let $U$ be an $r\times r$ unitary-valued matrix function such that 
{\em$\|H_U\|_{\text e}<1$}. The following are equivalent:

\noindent
{\em(i)} $U$ is very badly approximable;

\noindent
{\em(ii)} the Toeplitz operator $T_{zU}:H^2(\C^r)\to H^2(\C^r)$
has dense range in $H^2(\C^r)$;

\noindent
{\em(iii)} the Toeplitz operator $T_{\bar zU^*}:H^2(\C^r)\to H^2(\C^r)$
has trivial kernel.
\end{thm}

\Pf First of all it is trivial that (ii) is equivalent to (iii). 
The implication (i)$\Rightarrow$(ii) is proved in the Remark after the
proof of Theorem \ref{t4.3}.

It remains to show that (ii) implies (i). Again, it is explained in the Remark following the 
proof of Theorem \ref{t4.3} that $T_U$ is Fredholm.
Consider a Wiener--Hopf factorization of $U$:
$$
U=\Psi_2^*\left(\begin{array}{cccc}z^{d_1}&0&\cdots&0\\0&z^{d_2}&\cdots&0\\
\vdots&\vdots&\ddots&\vdots\\0&0&\cdots&z^{d_r}\end{array}\right)\Psi_1,
$$  
where $\Psi_1^{\pm1},\Psi_2^{\pm1}\in H^2(\mm_{r,r})$. It is well known and it is
easy to see that (ii) is equivalent to the condition that all Wiener--Hopf indices
$d_j$ are negative. It follows from the results of \S 3 of [AP2] that the 
superoptimal singular values $t_0(U),\cdots,t_{r-1}(U)$ are equal to 1 which means
that $U$ is very badly approximable. $\bl$

\

\

\noindent
Department of Mathematics
\newline
Kansas State University
\newline
Manhattan, Kansas 66506
\newline
USA

\end{document}